\documentclass[oneside,american]{amsart}
\usepackage{times}
\usepackage[T1]{fontenc}
\usepackage[latin1]{inputenc}
\usepackage{a4wide}
\usepackage{setspace}
\usepackage{amssymb}

\usepackage{psfrag}
\usepackage{float}
\usepackage{graphicx}
\usepackage{bbm} 
\usepackage{color}

\makeatletter
 \theoremstyle{plain}    
 \newtheorem{thm}{Theorem}[section]
 \numberwithin{equation}{section} 
 \numberwithin{figure}{section} 
 \theoremstyle{plain}
 \theoremstyle{plain}    
 \newtheorem*{thm*}{Theorem} 
 \newtheorem*{mthm*}{Main Theorem} 
 \theoremstyle{plain}    
 \newtheorem{lem}[thm]{Lemma} 
 \theoremstyle{remark}
 \newtheorem*{rem*}{Remark}
 \theoremstyle{definition}
 \newtheorem{defn}[thm]{Definition}
 \theoremstyle{plain}    
 \newtheorem{prop}[thm]{Proposition} 
 \theoremstyle{plain}    
 \newtheorem{cor}[thm]{Corollary} 

\DeclareMathOperator{\sign}{sign}
\DeclareMathOperator{\Lim}{Lim}

\renewcommand{\Re}{\mathfrak{Re}}
\renewcommand{\Im}{\mathfrak{Im}}
\renewcommand{\tilde}{\widetilde}
\renewcommand{\hat}{\widehat}

\newcommand{\Z}{\mathbb{Z}}

\newcommand{\R}{\mathbb {R}}
\newcommand{\Q}{\mathbb {Q}}

\renewcommand{\H}{\mathbb {H}}
\newcommand{\C}{\mathbb {C}}

\newcommand{\N}{\mathbb {N}}
\newcommand{\K}{\mathbb {K}}

\newcommand{\I}{\mathbb {I}}

\AtBeginDocument{ 
}

\usepackage{babel}
\makeatother
\begin{document}

\date{\today{}}

\author{Marc Kesseb\"{o}hmer and Bernd O. Stratmann}

\address{Fachbereich 3 - Mathematik und Informatik, Universit\"{a}t Bremen,
D--28359 Bremen, Germany}

\email{mhk@math.uni-bremen.de}

\address{Mathematical Institute, University of St Andrews, St Andrews KY16
9SS, Scotland}

\email{bos@maths.st-and.ac.uk}

\title{
Limiting modular symbols and their fractal geometry}



\begin{abstract}
In this paper we use  fractal geometry to investigate  boundary 
aspects of the first homology group for finite coverings 
of the modular surface. We obtain a complete 
description of algebraically invisible parts of this homology group. 
More precisely, we first show that for any
modular subgroup the geodesic forward dynamic on
the associated  surface admits a canonical symbolic 
representation by a finitely irreducible shift space.  We then use this representation to derive an `almost complete' 
 multifractal description of the higher--dimensional level sets arising from
Manin--Marcolli's limiting modular symbols.  
\end{abstract}

\maketitle

\section{Introduction }
\let\languagename\relax
Let  $\mathcal{C}_{2}(G)$ refer to the space of cusp forms 
 of weight $2$ for some arbitrary modular subgroup $G$. 
 That is, $G$ is a finite index 
subgroup of the modular group $\Gamma:=\textrm{PSL}_{2}(\Z)$.
It is well known that there is a dual pairing between 
$\mathcal{C}_{2}(G)$ and the first homology group $H_{1}(M_{G},\R)$ of 
the associated compactified cusped
surface $M_{G}$ of genus 
$\mathfrak{g}$. That is, we have
\[
\langle\,\cdot\,,\,\cdot\,\rangle:H_{1}(M_{G},\R)\times\mathcal{C}_{2}(G)\rightarrow\C,\;\textrm{where }
\;\langle\gamma,f\rangle:=\int_{\gamma}f(z)\, dz.\]
Each element of $H_{1}(M_{G},\R)$ can be represented by integrating 
 the $1$--form $f \, dz$ along some  smooth 
path between two points $\xi, \eta$ in $\H \cup 
P^{1}\left(\Q\right)$, and this determines the  modular symbol $\{\xi, \eta\}_{G} \in H_{1}(M_{G},\R)$.  
A possible way to extend these symbols to the non--cuspital 
boundary of hyperbolic space, and therefore to give a non--trivial 
homological meaning to  algebraically  invisible parts  
 of  $H_{1}(M_{G},\R)$, has been suggested by Manin and Marcolli
 \cite{ManinMarcolli:01}. 
 They introduced the concept  `limiting modular 
symbol', which is given for $x \in \R$
 by (whenever the limit exists)
\[ \ell_{G}(x):= \lim_{t\to \infty} \frac{1}{t} \{i,x+i \exp(-t)\}_{G}\in  H_{1}(M_{G},\R). \]
Note that  the limit in the definition of $\ell_{G}$  exists if and only if it exists for each 
$1$--form $f \, dz$  with $f \in \mathcal{C}_{2}(G)$,  and hence it is 
sufficient to compute it for a fixed complex basis 
$\hat{f}_{1},\ldots,\hat{f}_{\mathfrak{g}}$ of 
$\mathcal{C}_{2}\left(G\right)$. \\
The aim of this paper is to use fractal geometry in order to
investigate the level 
sets which arise naturally 
from these limiting modular 
symbols. That is, for $\alpha \in \R^{2\mathfrak{g}}$ we consider
 \[
\mathcal{F}_{\alpha}:=\left\{ 
x\in\R:\left(\langle\ell_{G}(x),f_{1}\rangle, \ldots, \langle\ell_{G}(x),f_{2 
\mathfrak{g}}\rangle\right)=\alpha\right\} ,\]
where $ f_{2k-1}:=\Re(\hat{f}_{k})$ and $ f_{2k}:=\Im(\hat{f}_{k})$, for 
$k=1,\ldots,\mathfrak{g}$. \\
A first analysis of this type of level 
sets was given in \cite{ManinMarcolli:01} and \cite{Marcolli:03} for modular 
subgroups
which satisfy the there so called `Red-condition' (see \cite{ManinMarcolli:01}).
There it was 
shown that for these groups  $\frac{1}{t} \{i,x+i \exp(-t)\}_{G}$ 
converges weakly to 
zero with respect to the Lebesgue measure on the unit interval. 
Subsequently, this result 
was improved  in \cite{Marcolli:03} by
showing that  $\ell_{G}(x)$ is equal to zero  Lebesgue--almost 
everywhere.  Besides,  these papers obtained
``non-vanishing ''  of limiting modular symbols
\emph{only}  for  the end points of closed geodesics, that is  
for quadratic surds. In these trivial cases the limiting modular 
symbol turns out to be given by  integrating along the closed 
geodesic and then normalizing by the  hyperbolic length of 
that geodesic. \\  
The aim of this paper is to extend these results to arbitrary 
modular subgroups and to obtain  that the limiting 
modular symbol is not equal to zero for large classes of perfect sets of positive 
Hausdorff dimension. In particular, a side result of our analysis will be that 
every modular subgroup satisfies the Red-condition, and hence that 
the results of \cite{ManinMarcolli:01} and \cite{Marcolli:03} do in 
fact hold for arbitrary modular 
subgroups. \\
Our main results are summarized in the following theorem, where
$\hat{\beta}_{G}:\R^{2\mathfrak{g}} \rightarrow\R$
refers to the proper concave (negative) Legendre transform  of the
proper convex function $\beta_{G} :\R^{2\mathfrak{g}}\to\R$, given by 
$\hat{\beta}_{G}(\alpha):=\inf_{t\in\R^{2\mathfrak{g}}}\left(\beta_{G} \left(t\right)-
\left( \alpha|t\right) \right)$.

\begin{mthm*} For every modular subgroup $G$
 there exists a strictly convex, differentiable function
$\beta_{G}:\R^{2\mathfrak{g}}\rightarrow\R$ such that for each $\alpha\in\nabla\beta_{G}\left(
\R^{2\mathfrak{g}}\right) \subset\R^{2\mathfrak{g}}$,
\begin{eqnarray}\label{mainbeta}
\dim_{H}\left(\mathcal{F}_{\alpha}\right)=\hat{\beta}_{G}(\alpha).\end{eqnarray}
In here, we have that $\beta_{G}(0)=1$, and that $\beta_{G}$ has a 
unique minimum at $0$. Also, we 
 in particular have   \[
\ell_{G}(\mathcal{F}_{\alpha})=\left\{ h_{\alpha}\right\} ,\]
 where $h_{\alpha}\in H_{1}(M_{G},\R)$ is uniquely determined by
$\left(\langle h_{\alpha},f_{1}\rangle,\ldots,\langle h_{\alpha},f_{2\mathfrak{g}}\rangle\right)=\alpha$.
Furthermore, the description of the spectrum in 
(\ref{mainbeta}) is almost complete in the sense
that \[ 
\overline{\nabla\beta_{G}\left(\R^{2\mathfrak{g}}\right)} =
\left\{ 
\alpha\in\R^{2\mathfrak{g}}:\mathcal{F}(\alpha)\neq\emptyset\right\} . \]
\end{mthm*}
The paper is organized as follows. In Section 2 we recall some basic 
facts on modular symbols. This includes a brief histogram on 
exact dual pairings. In Section 3 we  first give a reminder on a 
 beautiful construction which allows to visualize real numbers from a 
modular surface perspective. Subsequently, we show how this can be 
 generalized to arbitrary modular subgroups. 
Here, the main result will be  to use ergodicity of the geodesic flow 
to obtain that the so obtained  generalized modular shift spaces are 
always necessarily finitely primitive.  In Section 4 we introduce  
limiting modular symbol naturally arising from these generalized 
modular shift spaces, and show how these relate to the underlying 
branched geometry of numbers.  In Section  5 we will collect 
facts from the thermodynamic formalism  which turn out to be crucial in the 
proof of our main theorem,  which  will then be given  in Section 6.

\begin{rem*}
    1.  It is well known that the modular subgroup quotient $\H / G$ can be 
    viewed as a complex algebraic curve  permitting an arithmetic 
    structure. Similar to the familiar picture in which $\H / 
    \Gamma $ represents the moduli space of elliptic curves, $\H / G$
    represents the moduli space $\mathcal{M}(G)$ of elliptic curves 
    $\mathcal{E}(G)$ equipped with some finite additional  structure 
    determined by $G \backslash \Gamma$. Hence, by adding the cusp 
    points we obtain a pre-compactification 
    of $\mathcal{M}(G)$, given by including all possible ways  
    in which $\mathcal{E}(G)$ degenerates to $\C \setminus \{0\}$. 
    Moreover, if we further  include  all degenerations of $\mathcal{E}(G)$
    to noncommutative tori,  we derive the compactified moduli space 
    $\mathcal{M}_{0}(G)$ whose boundary is given by the noncommutative 
    space $P^{1}(\R) / G$. Since the level sets 
    $\mathcal{F}_{\alpha}$ are clearly $G$--invariant, a transfer of 
    the results in this paper to the language of  elliptic curve 
    degenerations leads to that the modular multifractal spectrum 
    corresponds to a continuous family  of elements of the boundary of 
    $\mathcal{M}_{0}(G)$
    consisting of the `bad quotients' $\mathcal{F}_{\alpha}/G$.
    (For the relation of noncommutative  geometry and modular 
    subgroups, and for some of the literature on this, we refer to the recent 
    survey article \cite{CM:06}).
     
2.  Let us also mention that the concept `limiting modular symbol' could easily be 
extended to more general concepts of  `modular symbols at infinity'. 
For instance, one could consider $\ell_{G,\phi,\psi}$ given by
\[ \ell_{G,\phi,\psi} (x):= {\Lim}_{t\to \infty} \phi(x,t) \{ i, x+i 
\psi(x,t)\}_{G} ,\]
for functions $\phi, \psi :\R \times \R^{+}\to \R^{+}$ such that  
$\phi(x ,t)$ and $\psi(x,t)$ tend to zero for $t$ tending 
to infinity. Here, $\Lim$ represents either $\lim$,  
$\limsup$, or $\liminf$. However, in this paper we concentrate on
 Manin--Markollis's limiting modular symbols for which $ \phi(x,t) :=1/t$, 
$\psi(x,t):= \exp(-t)$
and ${\Lim} 
:=\lim$.

\end{rem*}

\let\languagename\relax

\section{Preliminaries for modular symbols}

For  a modular subgroup $G$ consider  the space $\mathcal{C}_{k}(G)$ of cusp forms $f:\mathbb{H}\rightarrow\C$
of weight $k\in\Z$, given by
\begin{itemize}
\item $f$ is holomorphic on $\mathbb{H}$ as well as in each cusp of $G$; 
\item $f=(g')^{k/2}\cdot(f\circ g)$ for all $g\in G$; 
\item $f$ vanishes at each cusp of $G$. 
\end{itemize}
Throughout,  let $M_{G}$ refer to
the  (possibly branched)  covering surface $(\mathbb{H}\cup P^{1}\left(\Q\right))/G$
of $M_{\Gamma}$ of genus $\mathfrak{g}$, assumed to be compactified by 
having added 
the cusps.\\
Recall that a $p$-chain is a formal sum $c_{p}=\sum_{i}k_{i}N_{i}$, 
where the $N_{i}$
are $p$-dimensional smooth oriented submanifolds of $M_{G}$, 
and the
coefficients $k_{i}$ are elements of some abelian group $\mathbb{K}$. The 
$p$--homology group $H_{p}(M_{G},\mathbb{K})$ of $M_{G}$ is then
given via cycles and boundaries by 
$H_{p}(M_{G},\mathbb{K}):=\left\{ c_{p}:\partial 
c_{p}=\emptyset\right\} /\left\{ \partial c_{p+1}\right\} $.
Obviously, we have that $H_{p}(M_{G},\mathbb{K})=0$,  for each $p>2$. In 
particular,  $H_{1}(M_{G},\Z)$ is obtained geometrically by taking all loops
in $M_{G}$ as generators and then factoring out the relation that
two loops are homologous, that is they  only differ  by some boundary. 
By triangulating
$M_{G}$ such that the directed edges of the triangulation represent
the generators of the abelian group $H_{1}(M_{G},\Z)$ (modulo the
relations given by zero-homologous edge cycles), each element of $H_{1}(M_{G},\Z)$
can be written as a $\Z$--linear combination of the directed edges.
Hence, $H_{1}(M_{G},\Z)$ is a finitely generated abelian group which is
always equal to either $\Z^{2\mathfrak{g}}$ or to a free product of 
$\Z^{2\mathfrak{g}}$ with some
torsion subgroups of the form $\Z_{2}$ and/or  $\Z_{3}$. Given that
 $H_{1}(M_{G},\Z)$ is known, one can then apply the
universal coefficient theorem to determine 
$H_{1}(M_{G},\mathbb{K})$. Indeed, we have that $H_{1}(M_{G},\mathbb{K})=
 H_{1}(M_{G},\Z)\otimes_{\Z}\mathbb{K}$,
and hence $H_{1}(M_{G},\mathbb{K})$ is  a free $\K$--module of dimension $2\mathfrak{g}$,
for any ring $\mathbb{K}$. In particular, this shows that the space $H_{1}(M_{G},\mathbb{R})$
is a real vector space of dimension $2\mathfrak{g}$. \\
In order to see how $H_{1}(M_{G},\mathbb{R})$ fits in with 
$\mathcal{C}_{2}\left(G\right)$, note that by de~Rham theory we have 
that the $1$--cohomology group $H^{1}(M_{G},\mathbb{R})$
is isomorphic to the de Rham cohomology
$H_{\textrm{dR}}^{1}\left(M_{G},\mathbb{R}\right):=\{$closed $1$-forms$
\} / \{$exact 
$1$-forms$\}$,
and hence defines a dual pairing
(see e.g. the survey article \cite{EGH80}) \[
\langle\,\cdot\,,\,\cdot\,\rangle:H_{1}(M_{G},\R)\times H_{\textrm{dR}}^{1}\left(M_{G},\mathbb{R}\right)\rightarrow\R,\;
\textrm{given by }\;\langle\gamma,\omega\rangle:=\int_{\gamma}\omega.\]
By Hodge decomposition, there exists an isometry between $H_{\textrm{dR}}^{1}
\left(M_{G},\mathbb{R}\right)$ and the space $ H_{\Delta}^{1}\left(M_{G}\right)$
of harmonic $1$--forms. By considering real and imaginary parts of
the pullbacks of these harmonic forms to $\H\cup 
P^{1}\left(\Q\right)$,
we finally obtain an isomorphic representation of $H_{\Delta}^{1}\left(M_{G}\right)$
by the holomorphic cusp forms $\mathcal{C}_{2}\left(G\right)$. 
It follows that $\mathcal{C}_{2}\left(G\right)$ is a $\mathfrak{g}$--dimensional
$\C$--vector space (and hence a $2\mathfrak{g}$--dimensional $\R$--vector
space) (cf. also  \cite{Shimura:71})).
Hence,  summarizing the above, there is an exact
dual pairing of homology and cusp forms, given by the $\R$--bi-linear
map 
\[
\langle\,\cdot\,,\,\cdot\,\rangle:H_{1}(M_{G},\R)\times\mathcal{C}_{2}(G)
\rightarrow\C,\;\textrm{given by }\;
\langle\gamma,f\rangle:=\int_{\gamma}f(z)\, dz.\]
Also, let us recall for later use that by de~Rham theory we have for a fixed basis $\left(\gamma_{k}\right)_{k=1}^{2\mathfrak{g}}$,
consisting of cycles in $H_{1}\left(M_{G},\mathbb{R}\right)$, that the set \[
\left\{ \left(\left\langle \gamma_{k},h_{1}\right\rangle ,\ldots,\left\langle \gamma_{k},h_{\mathfrak{g}}
\right\rangle \right):k=1,\ldots,2\mathfrak{g}\right\} \]
 is a basis for $\mathbb{C}^{\mathfrak{g}}$ over $\R$ (cf. \cite[p. 10]{Cremona:1992}).
For the $2\mathfrak{g}\times2\mathfrak{g}$ real period matrix $\Pi_{G}:=\left(\left\langle \gamma_{k},f_{j}\right\rangle \right)_{k,j}$,
this means that \begin{equation}
\Pi_{G}\in\textrm{GL}_{2\mathfrak{g}}\left(\mathbb{R}\right).\label{period}\end{equation}
Note that the above dual  pairing shows that an element of $H_{1}(M_{G},\mathbb{R})$
can be viewed as a path in $M_{G}$, or alternatively as a path in $\H\cup P^{1}\left(\Q\right)$.
Therefore, by viewing it as a path in $\H\cup P^{1}\left(\Q\right)$ 
and noting that only the end points $\xi,\eta 
\in \H\cup P^{1}\left(\Q\right)$ of the path matter, this allows to 
represent an elements of  $H_{1}(M_{G},\mathbb{R})$ by the so called modular
symbol $\left\{ \xi,\eta\right\} _{G}\in H_{1}(M_{G},\R)$. In 
particular, for each $g\in G$, any smooth path from $\xi\in\H\cup P^{1}\left(\Q\right)$
to $g(\xi)$ projects to a closed path in $M_{G}$, and hence corresponds
to a homology class in $H_{1}(M_{G},\Z)$. Clearly, this class is  represented
by the modular symbol $\{\xi,g(\xi)\}_{G}$, obtained
by integrating $1$--forms $f\, dz$, for $f\in\mathcal{C}_{2}(G)$,
along any smooth path from $\xi$ to $g(\xi)$. One easily verifies
that $\{\xi,g(\xi)\}_{G}$ does not depend on $\xi$, and that the
assignment $g\mapsto\{\xi,g(\xi)\}_{G}$ gives rise to a surjection of $G$
onto $H_{1}(M_{G},\Z)$ (the kernel of this group homomorphism is
generated by  the commutators
of $G$). For the calculus with modular symbols, the following immediate
identities are useful. For each $\xi,\eta,\zeta\in\mathbb{H}\cup P^{1}\left(\Q\right)$
and $g\in G$, we have \[
\{\xi,\xi\}_{G}=0,\{\xi,\eta\}_{G}=-\{\eta,\xi\}_{G},\{\xi,\eta\}_{G}+\{\eta,\zeta\}_{G}=
\{\xi,\zeta\}_{G},
\]
 \[
\{\xi,\eta\}_{G}=\{ g(\xi),g(\eta\}_{G},\;\textrm{and }\;\{\xi,g(\xi)\}_{G}=\{\eta,g(\eta)\}_{G}.\]
 So far, we only considered paths between points in $\mathbb{H}\cup P^{1}\left(\Q\right)$
which lie in a single $G$--orbit of some element of $\mathbb{H} \cup P^{1}\left(\Q\right)$.
In general, that is for arbitrary modular subgroups $G$, it is not clear how to define
modular symbols between elements of $P^{1}\left(\Q\right)$ which
are not in a single $G$--orbit. However, for congruence subgroups
$\Gamma_{0}(N)$, defined by \[
\Gamma_{0}(N):=\left\{ \left(\begin{array}{cc}
a & b\\
c & d\end{array}\right) \in\Gamma\,|\, c\equiv 0 \mod N \right\} 
\;\textrm{for }\; N\in\N ,\]
 it is well known that this can be resolved, and this is the essence
of the following theorem.
\begin{thm*}
[Manin--Drinfeld \cite{Manin:72,Drinfeld:73}] For each $\xi,\eta \in P^{1}\left(\Q\right)$,
we have \[
\{\xi,\eta\}_{\Gamma_{0}(N)}\in H_{1}\left(M_{\Gamma_{0}(N)},\Q\right).\]
\end{thm*}
Finally, let us recall a view useful facts about modular subgroups 
and in particular also  congruence subgroups
$\Gamma_{0}(N)$. \\
For instance, for the
index $\kappa_{N}:=[\Gamma:\Gamma_{0}(N)]$ of $\Gamma_{0}(N)$
in $\Gamma$, we have (\cite{Schoeneberg:74}) \[
\kappa_{N}=N\prod_{p|N}\left(1+\frac{1}{p}\right).\]
 Also, for the number $N_{k}$ of $\Gamma_{0}(N)$-inequivalent elliptic
fixed points of order $k\in\N$ and the number $N_{\infty}$ of $\Gamma_{0}(N)$-inequivalent
parabolic fixed points, we have (\cite{Shimura:71}) \[
N_{2}=\left\{ \begin{array}{ll}
0 & \:\textrm{if \,}4|N\\
\prod_{p|N}\left(1+\left[(-1):p\right]\right) & \:\textrm{else},
\end{array}\right.\quad N_{3}=\left\{ \begin{array}{ll}
0 & \textrm{if }\,9|N\\
\prod_{p|N}\left(1+\left[(-3):p\right]\right) & \textrm{else},\end{array}\right.\]
 and \[
N_{\infty}=\sum_{d|N,d>0}\phi\left(\textrm{g.c.d.}\left(d,N/d\right)\right),\]
 where $\phi$ is the Euler function and $[\,:\,]$ refers to the
Legendre symbol of quadratic residues. \\
For arbitrary 
modular subgroups $G$,  the following formula
for the genus $\mathfrak{g}$ of $M_{G}$
is an immediate
consequence of the Riemann-Roch Theorem. With $R_{k}$ the number
of $G$-inequivalent elliptic fixed points of order $k\in\N$, $R_{\infty}$
the number of $G$-inequivalent parabolic fixed points, and $\kappa$
the index of $G$ in $\Gamma$, we have  (\cite{Shimura:71}) \[
\mathfrak{g}=1+\frac{\kappa}{12}-\frac{R_{2}}{4}-\frac{R_{3}}{3}-\frac{R_{\infty}}{2}.\]

\section{Modular shift spaces}
One of the main ideas of this paper is to investigate limiting modular 
symbols by using a shift space which  generalizes  
the shift space for the continued fraction expansion of elements of
$[-1,1]$ (see e.g. \cite{Art:24, series:85, AF:91, GK:01}). 
In this section we give the construction of this shift space 
canonically associated with the geodesic dynamic on the Riemann 
surface arising from a modular subgroup.  
Note that this construction extends the usual coding procedure for the 
modular surface to arbitrary modular subgroups. For ease of notation,
we put $\mathcal{I}:=[-1,1] \cap \I, \mathcal{I}_{-1}:=
  [-1,0] \cap \I$ and $\mathcal{I}_{+1}:= [0,1]\cap \I$, where $\I$
  denotes the set of irrational numbers. 
\\
We begin with  recalling from \cite{series:85} the notion of
`type--change'  for a geodesic in the upper half--plane $\H$. 
For this note  that $\H$ can be tiled by the so called Farey 
tesselation, that is the tesselation by 
$\Gamma$--translates of the ideal triangle with cusp--vertices at $0, 
1$ and $\{i\infty\}$. Consequently, each oriented  geodesic $l$ 
with  irrational end points is covered by infinitely many tiles of 
this 
tesselation. Especially, by travelling along $l$ towards the positive 
direction each of these tiles gets intersected such that there is a 
single vertex of the three cusp-vertices always  seen either on the left or on the 
right of the intersection of $l$ with  the tile (the other 
two vertices are seen on the opposite side). In case the single vertex is seen  on the left, we say that 
the visit is of 
type $L$, otherwise it is called of type $R$. If in here two 
successive visits  are of different type, then one says that $l$ 
changes type at the point where it intersects the edge at which   
the two involved tiles intersect.\\
Now, let us consider the set $\tilde{\mathcal{L}}_{\Gamma}$ of oriented 
geodesics $l$ in $\H$ with 
initial point $l_{-}$ and end point 
$l_{+}$, given by
\[ \tilde{\mathcal{L}}_{\Gamma} := \{l=(l_{-},l_{+}): 
0<|l_{-}|\leq 1 \leq |l_{+}|, l_{-} \cdot l_{+}<0, \textrm{ and } \, l_{-} 
,l_{+} \in \I\}.\]
Each element 
$l$ of 
$\tilde{\mathcal{L}}_{\Gamma}$ can then be coded by its successive
`type--changes', that is
\[
l \, \,  \textrm{ is  coded by} \, \,   \left\{ \begin{array}{ll}
\ldots \ldots L^{n_{-2}} R^{n_{-1}} y _{l}L^{n_{1}}R^{n_{2}}\ldots 
\ldots &
\:\textrm{if \,} l_{+} \geq 1\\
\ldots \ldots R^{n_{-2}} L^{n_{-1}} y_{l} R^{n_{1}}L^{n_{2}} \ldots 
\ldots &
\:\textrm{if \,} l_{+} \leq -1,
\end{array}\right. \]
where $y_{l}$ refers to the point 
at which $l$ intersects the imaginary axis.  This type of coding 
is closely related to the continued fraction expansion 
 of elements $y \in  \mathcal{I}_{+1}$, given for 
$y_{1},y_{2},\ldots\in \N$ by
\[y=[y_{1},y_{2}, \ldots]:=\frac{1}{{\displaystyle y_{1}+\frac{1}{{\displaystyle 
 y_{2}+\;\ldots\;}}}}.\]
Namely,  we have that
\[ \begin{array}{ll}  
l_{-}=-[n_{-1}, n_{-2},\ldots ] \textrm{ \, and \, }  l_{+}=[n_{1}, 
n_{2}, \ldots 
]^{-1} &
 \:\textrm{if \,} l_{+} \geq 1  \\
 l_{-}= [n_{-1}, n_{-2},\ldots ] \textrm{ \, and \, }   
 l_{+}=-[n_{1}, n_{2}, \ldots 
]^{-1}  
& \: \textrm{if \,} l_{+} \leq -1. \end{array}\]
Next, consider the subset $\tilde{\mathcal{C}}_{\Gamma}$ of the unit 
tangent space $UT(\H)$ consisting of
all those unit tangent vectors which are based at the imaginary axis 
and which give rise 
to  geodesics $l \in 
\tilde{\mathcal{L}}_{\Gamma}$.
We then have that the 
Poincar\'e section $\tilde{\mathcal{S}}_{\Gamma}$ associated with 
$\tilde{\mathcal{L}}_{\Gamma}$  is given by  the
canonical projection of $\tilde{\mathcal{C}}_{\Gamma}$
onto $UT(M_{\Gamma})$. 
 More precisely, let  $l\in \tilde{\mathcal{L}}_{\Gamma}$ be given  such that $l$ is 
 coded by $\ldots L^{n_{-2}} R^{n_{-1}} y 
 _{l}L^{n_{1}}R^{n_{2}}\ldots$. With $T: z \mapsto z+1$ referring 
 to the parabolic generator of $\Gamma$, we have that $T^{-n_{1}}(l)$ 
 is a geodesic
  which starts in $[-(n_{1}+1),-n_{1}]$ and ends in $\mathcal{I}_{+1}$, and
  hence $T^{-n_{1}}(l)$
 is not an element of $\tilde{\mathcal{L}}_{\Gamma}$. However, if we additionally apply 
 the elliptic generator $S: z \mapsto -1/z$ of $\Gamma$,
 then we obtain that the resulting geodesic $l':=ST^{-n_{1}}(l)$  is 
 an element of $\tilde{\mathcal{L}}_{\Gamma}$, and one immediately verifies that 
 \[ l'_{-}=[n_{1}, n_{-1}, \ldots ] \textrm{ \, and \, }  
 l'_{+}=-[n_{2}, n_{3}, \ldots 
]^{-1}.\]
Hence, in this situation we have
\[
    ST^{-n_{1}}:  l=(-[n_{-1}, n_{-2}, \ldots ], 
 [n_{1}, n_{2}, \ldots 
] ) \mapsto l'=([n_{1}, n_{-1}, \ldots ], 
-[n_{2}, n_{3}, \ldots 
]).\]
The dynamical idea behind this coding step is as follows. Let 
$v_{l} \in \tilde{\mathcal{C}}_{\Gamma}$  
be given, and let $v'_{l}$ be the vector in 
 $\Gamma(\tilde{\mathcal{C}}_{\Gamma})$  obtained by sliding $v_{l}$ along 
 $l$ towards the 
 positive direction of $l$ until the  next type-change takes place.  The significance of 
 $ST^{-n_{1}}$ then is that  $ST^{-n_{1}}(v_{l}')$ is an element 
 of $\tilde{\mathcal{C}}_{\Gamma}$ such that its projection onto $UT(M_{\Gamma})$ is 
 precisely the first return  to $\tilde{\mathcal{S}}_{\Gamma}$ when starting
 from  the projection of $y_{l}$ onto $M_{\Gamma}$  towards the direction
 of $v_{l}$.  \\
 Clearly, we can proceed similar if $l$ turns out to be coded by 
$\ldots R^{n_{-2}} L^{n_{-1}} y_{l} R^{n_{1}}L^{n_{2}} \ldots$. In 
this case we obtain that
\[
    ST^{n_{1}}:  l=([n_{-1}, n_{-2}, \ldots ], 
 -[n_{1}, n_{2}, \ldots 
]^{-1} ) \mapsto l'=(-[n_{1}, n_{-1}, \ldots ], 
[n_{2}, n_{3}, \ldots 
]^{-1}).\]
This procedure is summarized by the  Poincar\'e--map 
$\tilde{\mathcal{P}}_{\Gamma}:\mathcal{L}_{\Gamma} \to 
\mathcal{L}_{\Gamma}$, 
given by
\[ \tilde{\mathcal{P}}_{\Gamma}(l):= 
\left\{ \begin{array}{ll}
 ST^{-n_{1}} (l) &
\:\textrm{if \,} l=(-[n_{-1}, n_{-2}, \ldots ], 
[n_{1}, n_{2}, \ldots 
]^{-1})\\
ST^{n_{1}} (l)  &
\:\textrm{if \,} l= ([n_{-1}, n_{-2}, \ldots ], 
-[n_{1}, n_{2}, \ldots 
]^{-1}).
\end{array}\right. \]
Here it is important to remark that  the restriction $\mathcal{P}_{\Gamma}$
of the action of $\tilde{\mathcal{P}}_{\Gamma}$
to the second coordinate can 
also be described by the `twisted Gauss--map' 
\[ \mathcal{G}_{\Gamma}:\mathcal{I} \to \mathcal{I},
x\mapsto S \mathcal{P}_{\Gamma}S(x).\]
The reason why  $\mathcal{G}_{\Gamma}$ is called twisted Gauss--map 
originates from its link to  the usual Gauss--map
$\mathcal{G}: \mathcal{I}_{+1} \to \mathcal{I}_{+1}, x \mapsto 1/x -
\left[\left[1/x\right]\right]$   (where $[[ 1/x]]$
denotes the integer part of $1/x$).
Namely, one immediately verifies 
\[\mathcal{G}_{\Gamma} ( :=S\mathcal{P}_{\Gamma}S) : x \mapsto -\sign(x) \mathcal{G}(|x|).\]
Note that $\mathcal{G}_{\Gamma}$ can alternatively be derived directly if
 instead of $\tilde{\mathcal{L}}_{\Gamma}$   one starts with the following 
 set of oriented $\H$--geodesics
 \[ \mathcal{L}_{\Gamma} := \{l=(l_{-},l_{+}): 0<|l_{+}|\leq 1 \leq |l_{-}|,
 l_{-} \cdot l_{+}<0, \textrm{ and } \, l_{-},l_{+} \in \I\}.\]
 These geodesics can then also be coded by the type--change mechanism  as 
 explained above. Here, the relevant section 
 $\mathcal{C}_{\Gamma} \subset UT(\H)$ is the set of unit tangent vectors 
 based at the imaginary axis, giving rise to the geodesics in 
 $\mathcal{L}$. The  Poincar\'e section arising from this will be 
 denoted by $\mathcal{S}_{\Gamma}$.  One then immediately verifies 
 that $\mathcal{G}_{\Gamma}$ coincides with the map obtained by 
 projecting the  Poincar\'e map associated with this alternative approach onto the 
 second coordinate. \\
Also, note  that  for
 $k \in \Z^{^{\times}}:=\Z\setminus\left\{                                                   
  0\right\}$  the $k$-th inverse branch 
$\mathcal{P}_{\Gamma,k}^{-1}$ of $\mathcal{P}_{\Gamma}$ has the property
\[ \mathcal{P}_{\Gamma,k}^{-1} : 
\left\{ \begin{array}{ll}
(-\infty, -1] \to \{y^{-1} \in [1,\infty): 
y=[k,y_{2},\ldots]\}&
\:\textrm{if \,} k \in \N \\
 \left[ 1,\infty \right) \to \{-y^{-1} \in (-\infty, -1]: 
y=[|k|,y_{2},\ldots]\} &
\:\textrm{if \,}  k \notin \N.
\end{array}\right. \]
In particular, this shows that $\mathcal{P}_{\Gamma,k}^{-1}$ can be 
expressed in terms of the generators of $\Gamma$ by
  $\mathcal{P}_{\Gamma,k}^{-1}= 
T^{k}S$.
On the other hand, for
the corresponding $k$-th inverse branch  $\mathcal{G}_{\Gamma,k}^{-1}$ 
of the twisted Gauss--map we immediately have 
\[\mathcal{G}_{\Gamma,k}^{-1}=  S\mathcal{P}_{\Gamma,k}^{-1}S
= ST^{k}SS= ST^{k}: 
\left\{ \begin{array}{ll}
\mathcal{I}_{+1} \to {I}_{-k} &
\:\textrm{if \,} k \in \N \\
 \mathcal{I}_{-1} \to {I}_{-k} &
\:\textrm{if \,}  k \notin \N,
 \end{array}\right. \]
where 
$I_{k}:=\{\sign (k)  [y_{1},y_{2},\ldots] \in \mathcal{I}:                                                 
  y_{1}=|k|\}$ refers to the basic intervals. \\
We can now use standard ergodic theory to obtain our actual code                                          
  space via the inverse branches of $\mathcal{G}_{\Gamma}$ as 
  follows (cf. \cite{Bowen:75, Ruelle:78}).  One immediately verifies that $\alpha                                                    
  :=\{I_{k}:k \in \Z^{^{\times}}\}$ is a partition of $\mathcal{I}$, such that                                         
  the sequence of refinements  $\left(\bigvee_{i=0}^{n-1}\mathcal{G}_{\Gamma}^{-i}(\alpha)                         
  \right)_{n\in \N}$
  generates                                              
  the Borel $\sigma$--algebra. Hence, in terms of inverse branches of 
  $\mathcal{G}_{\Gamma}$ the twisted continued fraction  coding of                                                 
  $\mathcal{I}$  is as follows. For $n_{1},n_{2},\ldots\in \N$, we have                                                          
  \[\begin{array}{ll}                                                                                              
  [n_{1}, n_{2},\ldots ]  = \lim_{k \to \infty}                                                                  
  ST^{-n_{1}}S T^{n_{2}}\ldots ST^{(-1)^{k}n_{k}} 
  \left(\mathcal{I}_{(-1)^{k}}\right)    \\                                
  -[n_{1}, n_{2},\ldots ]  = \lim_{k \to \infty}                                                                  
  ST^{n_{1}}S T^{-n_{2}}\ldots ST^{(-1)^{k+1}n_{k}} 
  \left(\mathcal{I}_{(-1)^{k+1}}\right) .                                        
\end{array}\]
Therefore, by defining  
the 
shift space 
\[ \Sigma_{*}:= \left\{(x_{1},x_{2},\ldots) \in \left(\Z^{^{\times}}\right)^{\N}: 
x_{i}x_{i+1}<0, \textrm{ for all } i \in \N \right\}\]
equipped with the shift map $ \sigma_{*}: (x_{1},x_{2},\ldots) \mapsto (x_{2}, 
x_{3},\ldots)$, 
 one 
immediately verifies 
that 
\[ \rho: \Sigma_{*} \to \mathcal{I},  \,
(x_{1},x_{2},\ldots) \mapsto  
-\sign(x_{1}) [|x_{1}|,|x_{2}|,\ldots] \] is a bijection 
for which  $\rho \circ \sigma_{*} = 
\mathcal{G}_{\Gamma} \circ \rho$. \\
Our next goal is to generalize this modular coding procedure to arbitrary 
modular 
subgroups $G$. For this,  let $E_{G}$ refer to a fixed set 
of left--coset representatives
of the quotient space $G\backslash\Gamma$. 
In 
this more general setting the relevant set of oriented
geodesics is given by $\mathcal{L}_{G}:= \bigcup_{e\in E_{G}} 
e(\mathcal{L}_{\Gamma})$. Note that there is a 1-1-correspondence between $\mathcal{L}_{G}$  
and
the Poincar\'e section $\mathcal{S}_{G}$  for the geodesic flow on $M_{G}$, 
where $\mathcal{S}_{G}$ is given by the canonical projection of $\mathcal{C}_{G}:=  \bigcup_{e\in 
E_{G}} e(\mathcal{C}_{\Gamma})$  onto $UT(M_{G})$. 
We then adopt the above modular coding 
procedure in order to obtain a code space also in this
more general situation. For this we proceed as follows.
Assume that $\overline{\Sigma}_{G} := \bigcup_{e\in E_{G}} \left(e(\mathcal{I}) 
\times \{e\}\right)$ is equipped with the topology inherited from $\R$.
The $G$--generalized twisted Gauss--map $\mathcal{G}_{G}:\overline{\Sigma}_{G}\to 
\overline{\Sigma}_{G}$ is then given by
\[ \mathcal{G}_{G}(x, e ):=  \left(eS \mathcal{P}_{\Gamma}Se^{-1}(x), 
e\right), \, 
\textrm{ for } \, e \in E_{G}, x \in e(\mathcal{I}) .\]
Analogous to the situation before, we now have that  for
 $k \in \Z^{^{\times}}$ and $e\in E_{G}$    the $(k,e)$-th inverse branch 
 $\mathcal{G}_{G, (k,e)}^{-1}$ of
 $\mathcal{G}_{G}$
 is given by 
 \[ \mathcal{G}_{G, (k,e)}^{-1} : \left\{\begin{array}{ll}
 e(\mathcal{I}_{+1})  \times \{e\} \to e(\mathcal{I}_{-1}) \times \{e\}, (x,e) 
 \mapsto \left(  eST^{k}e^{-1}(x),e\right) &
 \:\textrm{if \,}  k \in \N\\
  e(\mathcal{I}_{-1}) \times \{e\} \to e(\mathcal{I}_{+1}) \times \{e\}, (x,e) 
  \mapsto \left(eST^{k}e^{-1} (x), e \right) &
 \:\textrm{if \,} k \notin \N.
 \end{array} \right .\] 
Hence, we can again use standard ergodic theory to obtain our
actual code                                          
  space via the inverse branches of $\mathcal{G}_{G}$. 
   This time  the basic intervals 
  are 
  $I_{k,e}:=e(I_{k}) \times \{e\}$, for $k\in \Z^{^{\times}}$ and $e\in E_{G}$.
  Also, $\alpha_{G}                                                    
  :=\{I_{k,e}:k \in \Z^{^{\times}}, e \in E_{G}\}$ is a partition of 
  $\overline{\Sigma}_{G}$ such that                                         
  the sequence of refinements  $\left(\bigvee_{i=0}^{n-1}\mathcal{G}_{G}^{-i}(\alpha)                         
  \right)_{n\in \N}$  generates                                              
  the Borel $\sigma$--algebra of $\overline{\Sigma}_{G}$. Hence, in terms of inverse branches of 
  $\mathcal{G}_{G}$ the $G$--generalized twisted continued fraction  coding of                                                 
  $\overline{\Sigma}_{G}$  is as follows. Let $(x,e) \in \overline{\Sigma}_{G}$  be given. 
  Then $x \in e(\mathcal{I})$, and we have that there exist  $n_{1},n_{2},\ldots\in
  \N$   
  such that $e^{-1}(x) = \pm [n_{1},n_{2},\ldots]$.  By the above 
  modular coding, we then 
  have                                                          
  \begin{eqnarray}\label{approx} x =  \left\{ \begin{array}{ll}                                                                                              
 \lim_{k \to \infty}                                                                  
  e ST^{-n_{1}}S T^{n_{2}}\ldots ST^{(-1)^{k}n_{k}} 
  \left(\mathcal{I}_{(-1)^{k}}\right)   &
 \:\textrm{if \,} e^{-1}(x) \in \mathcal{I}_{+1}  \\                                
 \lim_{k \to \infty}                                                                  
  eST^{n_{1}}S T^{-n_{2}}\ldots ST^{(-1)^{k+1}n_{k}} 
  \left(\mathcal{I}_{(-1)^{k}}\right)   &
 \:\textrm{if \,} e^{-1}(x) \in \mathcal{I}_{-1} .                                       
\end{array} \right. 
\end{eqnarray} 
Clearly,  the assignment $
   \left(e(\pm [n_{1},n_{2},\ldots]) , e\right) \mapsto \left((\mp 
   n_{1}, \pm n_{2}, \ldots), e\right) $
   gives rise to a bijection  
   between  $\overline{\Sigma}_{G}$ and  $\tilde{\Sigma}_{G}:=\Sigma_{*}\times E_{G}$.  
   Unfortunately,   the space $\tilde{\Sigma}_{G}$ is not a proper 
   shift space. However, this can be achieved by keeping track of the
   cosets 
   $G e ST^{\pm n_{1}}S T^{\mp n_{2}}\ldots ST^{\pm n_{k}} $ which
   are 
   visited during the approximation of $x$ given in (\ref{approx}).
  That is, we successively mark down as a second parameter
  the cosets 
in which those images of the directed imaginary axis lie on which the 
type--changes occur. More precisely, we define   
the shift space $\Sigma_{G}$  by  \\
$\Sigma_{G}:=\{((x_{1},e_{1}),(x_{2},e_{2}),\ldots)\in(\Z^{^{\times}}\times
E_{G})^{\N}: (x_{1},x_{2},\ldots) \in \Sigma_{*},$
\[ \qquad \qquad \qquad\qquad \qquad \qquad\qquad \qquad \qquad  
\, e_{k+1}=\tau_{x_{k}}(e_{k}) \mbox{ for all } k\in\N \},\]
 where the map $\tau_{x_{k}}:E_{G}\to E_{G}$ is defined by, with 
 $\equiv_{G}$ referring to equivalence mod $G$,
 \[
\tau_{x_{k}}(e_{k}):\equiv_{_{G}}e_{k}ST^{x_{k}}.\]
One immediately verifies that the assignment 
\[ \left((x_{1},x_{2}, \ldots), e\right) \mapsto \left( 
(x_{1},e),(x_{2},\tau_{x_{1}}(e)), (x_{3},\tau_{x_{2}}(\tau_{x_{1}}(e))),\ldots 
\right)\]
defines an isomorphism between $\tilde{\Sigma}_{G}$ 
and $\Sigma_{G}$, and hence $\Sigma_{G}$ is also isomorphic to 
$\overline{\Sigma}_{G}$. Of course, the advantage in using 
$\Sigma_{G}$ to code the geodesic rays in $M_{G}$ which arise from 
$\mathcal{S}_{G}$  is that it becomes a proper shift space when  
 equipped  with the shift map \[ \sigma:\Sigma_{G}\rightarrow\Sigma_{G},
((x_{1},e_{1}),(x_{2},e_{2}),\ldots)\mapsto((x_{2},e_{2}),(x_{3},e_{3}),\ldots),\]
 as well as with the canonical metric $d$, given for $((x_{k},e_{k}))_{k},((x_{k}',e_{k}'))_{k}\in\Sigma_{G}$
by \[
d(((x_{k},e_{k}))_{k},((x_{k}',e_{k}'))_{k}):=\sum_{i=1}^{\infty}2^{-i}\,\left(1-\delta_{(x_{i},e_{i}),(x_{i}',e_{i}')}\right).\]
Note  that the system $(\Sigma_{G},\sigma)$ relates to ordinary continued
fraction expansions as follows. For $\left((x_{k},e_{k})\right)_{k}\in\Sigma_{G}$
one immediately verifies by way of finite induction, using the 
matrix representation of the elements in $\Gamma$, \begin{eqnarray*}
e_{k+1} & = & \tau_{x_{k}}(e_{k})\equiv_{_{G}}e_{k}ST^{x_{k}}=
\tau_{x_{k-1}}(e_{k-1})ST^{x_{k}}\\
 & \equiv_{_{G}} & e_{k-1}ST^{x_{k-1}}ST^{x_{k}} \equiv_{_{G}}\ldots\\
 & \equiv_{_{G}} & e_{1}ST^{x_{1}}\ldots ST^{x_{k}}=e_{1}
  \left(\begin{array}{cc}
 -\sign\left(x_{1}\right)p_{k-1}(x) 
 &  \left(-1\right)^{k}p_{k}(x) \\
 q_{k-1}(x) & \left(-1\right)^{k+1} \sign\left(x_{1}\right) 
  q_{k}(x) \end{array}\right).\end{eqnarray*}
 Here, $p_{n}(x)/q_{n}(x)$ refers to the $n$-th approximant of the ordinary 
 continued fraction expansion of $x:=[|x_{1}|,|x_{2}|,\ldots]$,
with the usual convention  $q_{0}(x) :=p_{-1}(x):=1$ and $q_{-1}(x):=p_{0}(x):=0$.

\begin{rem*}
1.   At  first sight it might appear that the step from  $\tilde{\Sigma}_{G}$ 
and/or $\overline{\Sigma}_{G}$
to $\Sigma_{G}$ is just technical and that it achieves only little.
However, this step will turn out to be crucial, since it will allow us to 
 employ  certain standard results from thermodynamic formalism, 
 a formalism which is well 
 elaborated for shift spaces of the type $\left(
 \Sigma_{G}, \sigma\right)$.
    
    2. We remark that $(\Sigma_{G},\sigma)$ can also be represented by the skew product 
$\left(\tilde{\Sigma}_{G}, \tilde{\sigma}\right)$, where 
$\tilde{\sigma}: \tilde{\Sigma}_{G}\rightarrow\tilde{\Sigma}_{G}$ is 
given by
\[
\tilde{\sigma}:
 \left((x_{1},x_{2},\ldots),e\right) \mapsto \left((x_{2},x_{3},\ldots),\tau_{x_{1}}(e)\right).\]
 One  immediately verifies that the assignment 
 $\tilde{\pi}((x_{1},e_{1}),(x_{2},e_{2}),\ldots):=((x_{1},x_{2},\ldots),e_{1})$
gives rise to an isomorphism $\tilde{\pi}:\Sigma_{G} 
\rightarrow  \tilde{\Sigma}_{G}$,
which is a dynamical conjugacy in the sense that $\tilde{\sigma}\circ \tilde{\pi}=\tilde{\pi} \circ\sigma$. \\
Also, note that throughout
we will often identify elements $((x_{k},e_{k}))_{k}\in\Sigma_{G}$ with 
elements $\left((x_{k})_{k},e_{1}\right)\in \tilde{\Sigma}_{G}$,
as well as with elements 
$\left(e_{1}(-\sign\left(x_{1}\right)\left[\left|x_{1}\right|,\left|x_{2}\right|,\ldots\right]),e_{1}\right)\in
\overline{\Sigma}_{G}$. 

3. Let us also mention, although we are not going to use
this in this paper, that $\left(\Sigma_{G}, \sigma\right)$ can  be 
represented by a conformal graph directed Markov system 
(for an extensive 
discussion of these systems, we refer to \cite{MauldinUrbanski:03}).
Namely, 
define
$\mathcal{V}:=\left\{ \left(e,\pm 1 \right):e\in E_{G}\right\} $ to 
be the finite set of vertices, $\mathcal{E}:=\Z^{^{\times}}\times E_{G}$
 the countable infinite set of edges, and let two functions $i,t:\mathcal{E}\to\mathcal{V}$
be given by $i\left(\left(k,e\right)\right):=\left(e,-\sign\left(k\right)\right)$
and 
$t\left(\left(k,e\right)\right):=\left(\tau_{k}\left(e\right),\sign\left(k\right)\right)$.
Furthermore, let the edge incidence matrix 
$A=\left(A_{u,v}\right)_{u,v \in \mathcal{E}}$  be defined by 
$A_{u,v}=1$ if 
$t\left(u\right)=i\left(v\right)$, and  $A_{u,v}=0$ otherwise. 
We then have that $\left(\mathcal{V},\mathcal{E},i,t,A\right)$
is a directed multigraph  with associated
incidence matrix $A$, and one immediately verifies that 
the subshift $\Lambda_{G}:=\left\{ \left(u_{k}\right)_{k} 
\in\mathcal{E}^{\N}:A_{u_{k},u_{k+1}}=1, \textrm{for all } k\in \N \right\} $
is isomorphic to $\Sigma_{G}$. In order to derive the conformal graph directed Markov 
 system,  define compact sets 
 $\mathcal{I}_{\left(e,+1\right)}:=e\left(\mathcal{I}_{+1}\right)$
and $\mathcal{I}_{\left(e,-1\right)}:=e\left(\mathcal{I}_{-1}\right)$, for
all $e\in E_{G}$. Also, for each $(k,e) \in \mathcal{E}$ define
\[
\phi_{\left(k,e\right)}:=e ST^{k}
\left(\tau_{k}\left(e\right)\right)^{-1}
:
\mathcal{I}_{\left(\tau_{k}\left(e\right),\sign\left(k\right)\right)} 
\to \mathcal{I}_{\left(e,-\sign\left(k\right)\right)},\]
where  we assume that the maps $\phi_{\left(k,e\right)}$ are represented as   
M\"obius transformations contained in  $\Gamma$. With these 
preparations, one  now immediately verifies that the system 
\[
\Phi_{G}:=\left\{ \phi_{u}:\mathcal{I}_{t\left(u\right)}\to 
\mathcal{I}_{i\left(u\right)} \, | \, u\in\mathcal{E}\right\} \]
satisfies all the requirements of a conformal graph directed
Markov system, apart from that a priori the maps $\phi_{u}$
are  not necessarily uniformly contracting. However, similar as
for Schottky groups (\cite[Example 5.1.5]{MauldinUrbanski:03}),
this can get resolved  by replacing the system $\Phi_{G}$ by
a sufficiently high iterate of 
itself.
\end{rem*}

The final aim of this section is to show that for any 
modular subgroup we have that the modular shift 
space $\Sigma_{G}$ satisfies a certain transitivity condition  called
`finitely irreducible' 
(for the definition, see Proposition \ref{irred} below). 
Let us remark that  the main results in \cite{ManinMarcolli:02} are based
on the assumption that  the  'Red--condition' holds (cf. \cite{ManinMarcolli:02}).
One immediately verifies that this
Red-condition is in fact equivalent 
to finite irreducibility of the shift space $\Sigma_{G}$.
  However, the  approach in \cite{ManinMarcolli:02}  allows to verify this 
  condition  only for congruence subgroups.
  The proof there (\cite[Proposition 1.2.1]{ManinMarcolli:02})
  is based on the fact that for congruence subgroups $\Gamma_{0}(N)$
  there is an isomorphism between the set of modular symbols and  
the set of M--symbols (that is $P^{1}(\Z/N \Z)$, the projective   
line 
over the ring of integers mod $N$ (see also \cite{Cremona:1992})).
This then allows to verify the Red--condition algebraically in terms of
elementary congruence calculations. In contrast to this, our approach  
is completely different. To obtain the result for \emph{all} modular 
subgroups $G$, we combine an elementary observation
for the shift space $\Sigma_{G}$ with the ergodicity of the geodesic flow
on $M_{G}$. \\
In the following,
$\Sigma_{G}^{n}$ refers to the set of admissible words of length $n$ in
the alphabet $\Z^{^{\times}}\times E_{G}$, and
$\Sigma_{G}^{\ast}:=\bigcup_{n\in\N}\Sigma_{G}^{n}$.

\begin{prop}\label{irred}
    For each modular subgroup $G$ we have that the modular shift space $\left(\Sigma_{G},\sigma\right)$
is finitely irreducible in the sense of \cite{MauldinUrbanski:03}).
That is, there exist a finite set $W\subset\Sigma_{G}^{\ast}$
such that for all $a,b\in\Z^{^{\times}} \times E_{G}$ there exist
$w\in W$ such that $awb\in\Sigma_{G}^{\ast}$. 
\end{prop}
\begin{proof}
 Let $(m,e'),(n,e'')\in\Z^{^{\times}}\times E_{G}$ be given. For 
 simplicity, let us only consider the case in which $m<0$ and $n>0$.
 Clearly, the remaining cases can be dealt with in an analogous way.
 Now, the aim is  to 
 show that there exists
       $c=c(G) \in \N$ such that $(m,e')w(n,e'')\in\Sigma_{G}^{\ast}$
       for some $w\in\Sigma_{G}^{j}$ with $j \leq c$. For this, 
       observe that with $\hat{e} := \tau_{m}(e')$ we have  
       that 
       $(m,e')(r,\hat{e})\in \Sigma_{G}^{2}$, for {\em all} $r\in \N$. 
       Likewise, observe that 
       if $(m,e')w(n,e'')\in\Sigma_{G}^{\ast}$ then 
       $(m,e')w(s,e'')\in\Sigma_{G}^{\ast}$, for {\em all} $s\in \N$.
       Combining these two observations, it follows that in order to 
       prove the assertion it is sufficient to show that 
       \begin{eqnarray}\label{assertion}  
	 \mbox{there exists  } w' \in \Sigma_{G}^{j-1}\mbox{ with }
       (r,\hat{e})w' (s,e'')\in\Sigma_{G}^{\ast},\mbox{ for {\em some} } r,s \in \N.
       \end{eqnarray}
For this,   note that we have by construction that  
$(k,e)\in\Z^{^{\times}}\times E_{G}$  represents the 
basic interval $I_{-k,e}$. Furthermore, in terms of cross 
sections we have 
that  $I_{-k,e}$  represents  a certain subset 
$\mathcal{C}_{k,e}$ of the cross section $e(\mathcal{C}_{\Gamma}) 
\subset \mathcal{C}_{G}$. That is, $\mathcal{C}_{k,e}$ is the set 
of  those unit--vectors $v$ which are 
based at $e(\{z \in \H: \Re(z)=0\}$
such that the oriented $\H$--geodesic  given by $v$
terminates in 
$e(I_{-k})$ and starts in either $e([1, 
\infty))$ (if $k$ is positive) or 
$e((-\infty , -1])$ (if $k$ is negative).
Define $\mathcal{C}_{e}^{-}:= \bigcup_{k\in  \N}  
\mathcal{C}_{k,e}$, and let $\mathcal{S}_{e}^{-} \subset \mathcal{S}_{G}$ 
be the projection 
of $\mathcal{C}_{e}^{-}$ onto $UT(M_{G})$. Expressing the assertion in 
(\ref{assertion}) in these terms,  it follows that we have to show that  
 there exists $v \in  
\mathcal{S}_{\hat{e}}^{-}$ and $v' \in  \mathcal{S}_{e''}^{-}$
such that $v'$ can be obtained by sliding $v$ in  positive 
direction along the $M_{G}$--geodesic given by $v$. But  this assertion 
follows immediately from
the fact that the geodesic flow on $M_{G}$ is ergodic.
This finishes the proof.
\end{proof}

\section{The limiting modular symbol for $\Sigma_{G}$}
 We already introduced  the  limiting modular 
symbol
$\ell_{G}$ in the introduction,
which there was  defined on $\R$. We now   define 
a slightly different version of such a symbol, namely the limiting
modular symbol $\tilde{\ell}_{G}$ which will be 
defined on $\Sigma_{G}$.
\begin{defn}
The limiting modular symbol $\tilde{\ell}_{G} : \Sigma_{G} \to
H_{1}(M_{G},\R)$ is defined for 
arbitrary $\left((x_{k},e_{k})\right)_{k}\in\Sigma_{G}$ by (whenever the limit exists
as an element of $H_{1}(M_{G},\R)$)
\[
\tilde{\ell}_{G}\left(\left((x_{k},e_{k})\right)_{k}\right):=
\lim_{t\rightarrow\infty}\frac{1}{t}\{ i,e_{1}(x+i\exp(-t))\}_{G}.\]
  Here, 
 we have set
 $x:=-\sign\left(x_{1}\right)\left[\left|x_{1}\right|,\left|x_{2}\right|,\ldots\right] \in \mathcal{I}$.
\end{defn}
Note that  $\ell_{G}(x)$ does not depend  on  the starting point $i$ 
of the paths along one integrates, nor does it depend on the choice of 
the geodesic $\{x+i \exp(-t): t \in \R \}$ (in fact, any path 
having $x$ as its only accumulation point in $P^{1}(\R)$ would do). 
Also, concerning the existence of the limit  in the definition of 
$\tilde{\ell}_{G}$ one can make the same remark as we made for $\ell_{G}$. 
That is, since $\langle\cdot,\cdot\rangle$ is a perfect dual pairing,
the existence of the limit  is guaranteed if  the following limit 
exists for each $f\in\mathcal{C}_{2}(G)$, or equivalently for each
member of a basis of $\mathcal{C}_{2}(G)$, \[
\lim_{t\rightarrow\infty}\left\langle \frac{1}{t}\{ i,e_{1}(x+i\exp(-t))\}_{G},f\right\rangle .\]
The following proposition gives the main result of this section. 
\begin{prop}\label{partial}
For $((x_{k},e_{k}))_{k}\in\Sigma_{G}$
we have \[
\tilde{\ell}_{G}\left(\left((x_{k},e_{k})\right)_{k}\right)=\lim_{n\rightarrow\infty}\frac{1}{2 \log 
q_{n}(|x|)}\sum_{k=1}^{n}
\{ e_{k}\left(i\infty\right),e_{k}\left(0\right)\}_{G}.\]
\end{prop}
\begin{proof}
\begin{figure}
\psfrag{el}{$e_1(l(x))$}\psfrag{ex}{$e_1 (x)$}
\psfrag{x1}{$\xi_1=y_1$}
\psfrag{x2}{$\xi_2$}\psfrag{x3}{$\xi_3$}
\psfrag{x4}{$\xi_4$}\psfrag{y2}{$y_2$}
\psfrag{y3}{$y_3$}\psfrag{y4}{$y_4$}
\psfrag{w1}{$\omega_1$}
\psfrag{w2}{$\omega_2$}\psfrag{w3}{$\omega_3$}
\includegraphics[%
  width=0.80\textwidth,
  keepaspectratio]{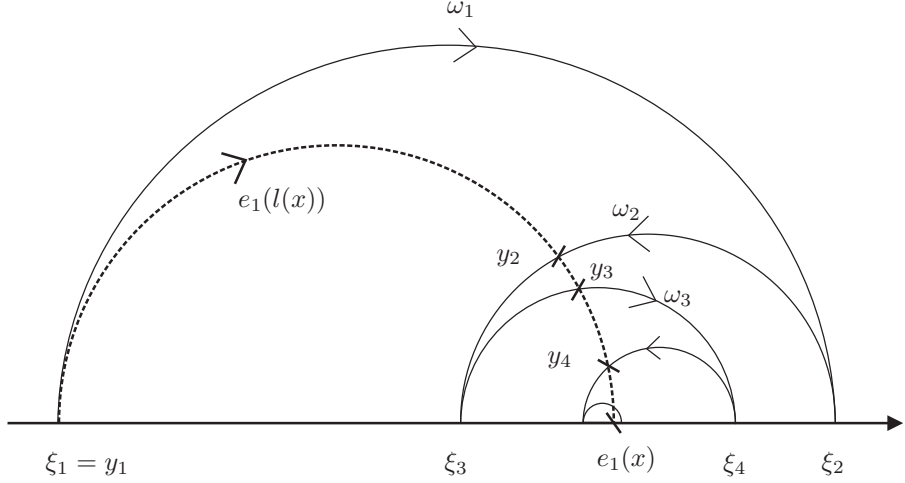}

\caption{Approximating path for a point $(x,e_1)\in \Sigma_G$.}\label{fig}
\end{figure}

Let $\left((x_{k},e_{k})\right)_{k}\in\Sigma_{G}$ be given. Our first aim 
is to show that   
\[\tilde{\ell}_{G,q}\left(\left((x_{k},e_{k})\right)_{k}\right):=\lim_{n\rightarrow\infty}
\frac{1}{2\log q_{n}(|x|)}\sum_{k=1}^{n}
\left\{ e_{k}\left(i\infty\right),e_{k}\left(0\right)\right\} _{G}\]
exists if and only if  $\lim_{n \rightarrow\infty} 1/t_{n} \{ 
i,e_{1}(x+i\exp(-t_{n}))\}_{G}$ exists, for some sequence $\left(t_{n}\right)_{n \in \N}$ 
 tending to infinity. More precisely, we will show 
that if one of these limits exists then both limits coincide, that is 
\begin{eqnarray}\label{qn}  \tilde{\ell}_{G,q} \left(\left((x_{k},e_{k})\right)_{k}\right) =
    \lim_{n \rightarrow\infty}\frac{1}{t_{n}}\{ 
i,e_{1}(x+i\exp(-t_{n}))\}_{G}.
\end{eqnarray}
For this we proceed similar to \cite[Proof of 
Theorem 0.2.1]{ManinMarcolli:02} as follows.
Let $l(x)$ refer to the 
oriented hyperbolic geodesic from $i \infty $ to $x$,  and define 
 $\xi_1=e_1(i\infty)$ and $\xi_{n}:= e_{1}\left(-\sign(x_{1})  p_{n-2}(|x|)/q_{n-2}(|x|)\right)$, 
 for 
$n \geq 2$. 
Then consider the path $\omega:= \omega_{1}\omega_{2} \ldots$ 
which runs in succession through the oriented hyperbolic geodesics
$\omega_{n}$ which start at  $\xi_{n}$ and end in 
$\xi_{n+1}$ (cf. Fig. \ref{fig}).  Clearly, when viewing $\omega$ as a path in $\H 
\cup P^{1}(\Q)$ 
it is a connected oriented path which approximates $e_{1}(x)$ in 
its forward direction. Next, define $y_{n}:= \omega_{n} \cap 
e_{1}(l(x))$ for $n \in \N$, and observe that the oriented geodesic path from $y_{n}$ 
to $y_{n+1}$ is  homologous to the geodesic path which runs from  $y_{n}$ 
via $\xi_{n+1}$ to $y_{n+1}$. It follows
\[  \left\{ y_{n}, y_{n+1} \right\}_{G} = \left\{y_{n}, \xi_{n+1} 
\right\}_{G} + \left\{\xi_{n+1},y_{n+1} \right\}_{G}, \, \textrm{ for 
all } n \in \N.\]
Before we continue with this argument, first observe that we have for all $n \in \N$,
\[ \left\{ \xi_{n}, \xi_{n+1} \right\}_{G}  
=                                                                                                 
  \left\{ e_{n}(i \infty), e_{n}(0) \right\}_{G} \, \textrm{ and }  
  \, 
  e_{n}(0) = e_{n+1}(i\infty). 
  \]
 Indeed, this can be seen as follows. Define $g_{1}:= \mathrm{id.}$, and 
 for $n \in \N$ let
 \[ g_{n+1}:=ST^{x_{1}}\ldots ST^{x_{n}}=\left(\begin{array}{cc}
 -\sign\left(x_{1}\right)p_{n-1}(|x|) & \left(-1\right)^{n}p_{n}(|x|)\\
 q_{n-1}(|x|) & 
 \left(-1\right)^{n+1}\sign\left(x_{1}\right)q_{n}(|x|)\end{array}\right).\]
  We then have that $e_{1} g_{n}(i\infty)= 
  \xi_{n} $ and
  $e_{1}g_{n}(0)=\xi_{n+1} $.
 Also, since $e_{n}\equiv_{G} e_{1}g_{n}$, there exists $\tilde{g}_{n}\in G$ such
 that $\tilde{g}_{n}e_{n}=e_{1}g_{n}$. Using these facts as well as the $G$--invariance
 of the modular symbol, we obtain
 \begin{eqnarray*}
 \left\{ e_{n}(i\infty),e_{n}(0)\right\} _{G}
  &=& \left\{ 
  \tilde{g}_{n}e_{n}(i\infty),\tilde{g}_{n}e_{n}(0)\right\} _{G} =
  \left\{ e_{1}g_{n}(i\infty),e_{1}g_{n}(0)\right\} _{G} \\
 & = & \left\{\xi_{n},\xi_{n+1}\right\}_{G}.
\end{eqnarray*}
Using this observation, we proceed with the above argument as 
follows. For  each $n \in \N$, we have
\begin{eqnarray*}
    \left\{ i, y_{n+1} \right\}_{G}  & = &     \left\{ i, y_{2} 
    \right\}_{G} + \left\{ y_{2}, y_{n+1} \right\}_{G} = \left\{ i, 
    y_{2} \right\}_{G} + \sum_{k=2}^{n} \left\{ y_{k}, y_{k+1} 
    \right\}_{G} \\
    & = & \left\{ i,  
    y_{2} \right\}_{G} + \sum_{k=2}^{n} \left( \left\{ y_{k}, \xi_{k+1} 
    \right\}_{G}  + \left\{ \xi_{k+1}, y_{k+1} 
    \right\}_{G} \right) \\
    &=& \left\{ i, y_{2}\right\}_{G} - \left\{ \xi_{2}, 
    y_{2} \right\}_{G}  -\left\{ y_{n+1}, 
    \xi_{n+2} \right\}_{G} + \sum_{k=2}^{n+1}  \left\{ \xi_{k}, 
    \xi_{k+1} 
    \right\}_{G}  \\
    & = & \left\{ i, 
    \xi_{1} \right\}_{G}  -\left\{ y_{n+1}, 
    \xi_{n+2} \right\}_{G}  
    +  \sum_{k=1}^{n+1}  \left\{ \xi_{k}, 
    \xi_{k+1} 
    \right\}_{G}  \\
    & = & \left\{ i, 
    \xi_{1} \right\}_{G}  -\left\{ y_{n+1}, 
    \xi_{n+2} \right\}_{G} 
	+  \sum_{k=1}^{n+1}  \left\{e_{k}(i\infty),e_{k}(0)\right\}_{G}. 
    \end{eqnarray*}
    Now, let $t_{n}$ be defined implicitly by $e_{1}(x+i 
    \exp(-t_{n})) := y_{n}$. Using elementary hyperbolic geometry 
    in the context of for instance  Ford circles
    (or alternatively, see e.g. \cite[paragraph 3]{KesseboehmerStratmann:04c}),   
    one immediately verifies that for all $n \in  \N$ sufficiently large
    we have $\exp(t_{n}) \asymp \left(q_{n}(|x|)\right)^{2}$.

This allows to finish the proof of 
(\ref{qn}) as follows.
\begin{eqnarray*}
    \tilde{\ell}_{G,q} 
    \left(\left((x_{k},e_{k})\right)_{k}\right) &=& 
    \lim_{n\rightarrow\infty}
    \frac{1}{2\log q_{n}(|x|)}\sum_{k=1}^{n}\left\{
    e_{k}\left(i\infty\right),e_{k}\left(0\right)\right\} _{G} \\
    &=&  \lim_{n\rightarrow\infty} \frac{1}{t_{n}} \left(  \left\{ i, y_{n} \right\}_{G} 
 +\left\{ y_{n}, 
	\xi_{n+1} \right\}_{G}   -\left\{ i, 
	\xi_{1} \right\}_{G} \right) \\
	&=&  \lim_{n\rightarrow\infty} \frac{1}{t_{n}}  \left\{ i, y_{n} 
	\right\}_{G} = \lim_{n \rightarrow\infty}\frac{1}{t_{n}}\{ 
i,e_{1}(x+i\exp(-t_{n}))\}_{G}.
\end{eqnarray*}
In order to finish the proof of the proposition, it remains to show 
that $\lim_{n \rightarrow\infty}\frac{1}{t_{n}}\{ 
i,e_{1}(x+i\exp(-t_{n}))\}_{G}$ is independent of the particular 
chosen sequence $\left(t_{n}\right)$. That is, our final aim is to 
show that the existence of $\tilde{\ell}_{G,q} 
\left(\left((x_{k},e_{k})\right)_{k}\right) $  implies that 
$\tilde{\ell}_{G,q} 
\left(\left((x_{k},e_{k})\right)_{k}\right) = \tilde{\ell}_{G} 
\left(\left((x_{k},e_{k})\right)_{k}\right) $.   
In order to prove this,  we argue similar as in 
\cite[paragraph 3]{KesseboehmerStratmann:04c}  as follows.
 Suppose that $\tilde{\ell}_{G,q} 
\left(\left((x_{k},e_{k})\right)_{k}\right) $ exists, and define 
$n_{t}:=\sup\left\{ n\in\N: 2 \log q_{n}(|x|)\leq t\right\} $
and $\alpha_{h}:=\left\langle \tilde{\ell}_{G,q}
\left(\left((x_{k},e_{k})\right)_{k}\right),h\right\rangle $, for 
arbitrary $t>0$ and $h \in \mathcal{C}_{2}(G)$.
We then have\\
\\
${\displaystyle \limsup_{t\to\infty}\left|\frac{\left\langle \left\{ i,e_{1}(x+i\exp(-t))\right\} _{G},h
\right\rangle }{t}-\frac{\left\langle \sum_{k=1}^{n_{t}}\left\{ e_{k}(i\infty),e_{k}(0)\right\} _{G},h
\right\rangle }{2\log q_{n_{t}}(|x|)}\right|}$\begin{eqnarray*}
 & = & \limsup_{t\to\infty}\left|\frac{2 \log q_{n_{t}}(|x|)\left\langle \left\{ i,e_{1}(x+i\exp(-t))\right\} _{G},
 h\right\rangle -t\left\langle \sum_{k=1}^{n_{t}}\left\{ e_{k}(i\infty),e_{k}(0)\right\} _{G},h
 \right\rangle }{2 t\log q_{n_{t}}(|x|)}\right|\\
 & \leq & \limsup_{t\to\infty}\left|\frac{\left\langle \left\{ i,e_{1}(x+i\exp(-t))\right\} _{G}-\sum_{k=1}^{n_{t}}
 \left\{ e_{k}(i\infty),e_{k}(0)\right\} _{G},h\right\rangle }{t}\right|\\
 &  &  \qquad\qquad\qquad\qquad+\limsup_{t\to\infty}\left|
 \frac{2\log q_{n_{t}}(|x|)-t}{t}\right|\left|\frac{\left\langle \sum_{k=1}^{n_{t}}\left\{ e_{k}(i\infty),e_{k}(0)\right\} _{G},
 h\right\rangle }{2\log q_{n_{t}}(|x|)}\right|\\
 & \leq & \limsup_{t\to \infty} \frac{\textrm{const.}}{t} + 
  \limsup_{n\to\infty} \left|\alpha_{h}\right| 
 \frac{\log |x_{n+1}|}{\log q_{n}(|x|)} =   0+
  \limsup_{n\to\infty}\left|\alpha_{h}\right|  \frac{\log |x_{n+1}|}{ \log q_{n}(|x|)}.\end{eqnarray*}
This shows that  if  $\alpha_{h}=0$ holds for all $h\in\mathcal{C}_{2}\left(G\right)$,
then $\tilde{\ell}_{G}\left(\left((x_{k},e_{k})\right)_{k}\right)$ exists
and has to be equal to 
$\tilde{\ell}_{G,q}\left(\left((x_{k},e_{k})\right)_{k}\right)$. 
Hence, in this case the proof is finished. Therefore, we can now assume  
without loss of generality that there exists $h\in\mathcal{C}_{2}\left(G\right)$ such 
that  $\alpha_{h}>0$. By the above, in order to finish the proof of 
the proposition it is sufficient to show 
that  $\limsup_{n\to\infty}\frac{\log 
|x_{n+1}|}{\log q_{n}(|x|)}=0$.
For this, observe \begin{eqnarray*}
\alpha_{h} & = & \lim_{n\to\infty}\frac{\left\langle \sum_{k=1}^{n+1}\left\{ 
e_{k}(i\infty),e_{k}(0)\right\} _{G},h\right\rangle }{2\log q_{n+1}(|x|)}
 \\ & = & \lim_{n\to\infty}\frac{\left\langle \sum_{k=1}^{n}\left\{ 
 e_{k}(i\infty),e_{k}(0)\right\} _{G}+\left\{ e_{n+1}(i\infty),e_{n+1}(0)\right\} _{G},
 h\right\rangle }{2\log q_{n}(|x|)+2 \log |x_{n+1}|}\\
 & = & \lim_{n\to\infty}\frac{\left\langle \sum_{k=1}^{n}\left\{ e_{k}(i\infty),e_{k}(0)\right\} _{G},
 h\right\rangle \left(1+\frac{\left\langle 
 \left\{ e_{n+1}(i\infty),e_{n+1}(0)\right\} _{G},h\right\rangle }{\left\langle \sum_{k=1}^{n}\left\{ e_{k}(i\infty),e_{k}(0)\right\} _{G},h\right\rangle }
 \right)}{2 \log q_{n}(|x|)\left(1+\frac{\log |x_{n+1}|}{\log 
 q_{n}(|x|)}\right)}\\
 & = & \alpha_{h} \lim_{n\to\infty}\frac{1+\frac{\left\langle \left\{ 
 e_{n+1}(i\infty),e_{n+1}(0)\right\} _{G},h\right\rangle }{\left\langle 
 \sum_{k=1}^{n}\left\{ e_{k}(i\infty),e_{k}(0)\right\} _{G},h\right\rangle 
 }}{1+\frac{\log |x_{n+1}|}{\log q_{n}(|x|)}}.\end{eqnarray*}
Now, suppose by way of contradiction  that $\limsup_{n\to\infty}\frac{\log |x_{n+1}|}{\log q_{n}(|x|)}>0$.
Then there exists a subsequence $\left(n_{k}\right)$ such that 
$\lim_{k\to\infty}\frac{\log |x_{n_{k}+1}|}{\log q_{n_{k}}(|x|)} >0$,
and consequently we have $\lim_{k\to \infty} |x_{n_{k}+1}|=\infty$. Combining 
this with our assumption 
$\alpha_{h}> 0$, it follows  \[
1=\lim_{k\to\infty}\frac{\log q_{n_{k}}(|x|)}{\left\langle 
\sum_{m=1}^{n_{k}}\left\{ e_{m}(i\infty),e_{m}(0)
\right\} _{G},h\right\rangle }\frac{\left\langle \left\{ 
e_{n_{k}+1}(i\infty),e_{n_{k}+1}(0)\right\} _{G},h\right\rangle }{\log 
|x_{n_{k}+1}|}=\frac{1}{\alpha_{h}} \cdot 0=0.\]
This is a contradiction, and hence it follows that
$\limsup_{n\to\infty}\frac{\log |x_{n+1}|}{\log q_{n}(|x|)}=0$.
\end{proof}
For our final result in this section, recall from the introduction that 
$f_{1},\ldots, 
 f_{2 \mathfrak{g}}$ refers to a fixed $\R$--basis of 
 $\mathcal{C}_{2}\left(G\right)$ given by  the real and imaginary part 
 of some complex basis 
 of $\mathcal{C}_{2}\left(G\right)$. We then define for $e_{1} \in E_{G}$ 
 and $\alpha \in \R^{2 \mathfrak{g}}$,\\
 $ \tilde{\mathcal{F}}_{\alpha}(e_{1}):=\big\{ 
x\in \mathcal{I} : ((x_{k},e_{k}))_{k} \in \Sigma_{G}$
\[ \qquad\qquad\qquad \mbox{such that } \left.\, \left(\langle\tilde{\ell}_G(((x_{k},e_{k}))_{k}),f_{1}\rangle,\ldots ,
  \langle\tilde{\ell}_G(((x_{k},e_{k}))_{k}),f_{2\mathfrak{g}}\rangle
  \right)=\alpha \right\} ,\]
   we have set
   $x:=-\sign\left(x_{1}\right) \left[\left|x_{1}\right|,\left|x_{2}\right|,\ldots \right]$.
\begin{lem}\label{eql}
   For each $e,e'\in E_{G}$ and $ \alpha \in \R^{2\mathfrak{g}}$, we 
   have
    \[
    \dim_{H}\left(\tilde{\mathcal{F}}_{\alpha}(e)\right)=\dim_{H}\left(\tilde{\mathcal{F}}_{\alpha}(e')
    \right).\]
    \end{lem}
\begin{proof}
    Let $e,e' \in E_{G}$ be given. Since $\Sigma_{G}$ is
    finitely irreducible, it follows that there exists $n \in \N$ and 
    $\left((x_{1},e_{1}) , \ldots ,(x_{n},e_{n})\right)\in \Sigma_{G}^{*}$ 
    such that $e_{1}=e$ and $ e_{n}=e'$. 
    This implies that there exists $g\in G$ such that $eST^{x_{1}}\ldots ST^{x_{n}}=ge'$. 
    Then note that for $\tilde{g}:
    =e^{-1}ge'=ST^{x_{1}}\ldots ST^{x_{n}}$
    one immediately verifies that $\tilde{g}(\mathcal{I})\subset\mathcal{I}$. 
    Using this observation, the $G$--invariance of the modular 
    symbol, and 
    the fact that the limiting modular symbol does not depend on the 
    starting point of the path along which one integrates, we obtain 
    for each $y \in \tilde{\mathcal{F}}_{\alpha}(e')$ and $\alpha \in \R^{2\mathfrak{g}}$, \begin{eqnarray*}
     \lim_{t\rightarrow\infty}\frac{1}{t}\left\{ i,e'(y+i\exp(-t))\right\} 
     _{G}& =& \lim_{t\rightarrow\infty}\frac{1}{t}\left\{ g(i),ge'(y+i\exp(-t))\right\} _{G}\\
     & = & \lim_{t\rightarrow\infty}\frac{1}{t}\left\{ 
     i,ee^{-1}ge'(y+i\exp(-t))\right\} _{G}\\ &= & \lim_{t\rightarrow\infty}\frac{1}{t}
     \left\{ i,e\tilde{g}(y+i\exp(-t))\right\} _{G}\\
     & = & \lim_{t\rightarrow\infty}\frac{1}{t}\left\{ i,e(\tilde{g}(y)+i\exp(-t))\right\} _{G}.
    \end{eqnarray*}
    This shows that 
    $\tilde{g}\left(\tilde{\mathcal{F}}_{\alpha}(e')\right)\subset 
    \tilde{\mathcal{F}}_{\alpha}(e)$. 
    Since $\tilde{g}$ is conformal, and hence in particular
    bi--Lipschitz, and since $e,e'\in E_{G}$ were arbitrary, 
     the lemma follows.
    \end{proof}

\section{Modular potential and pressure function}\label{ModularPotential}

In this section we collect results from the general thermodynamic 
formalism which  will be required in the proof of our Main Theorem.\\
Let $I:\Sigma_{G}\rightarrow\R$ refer to the canonical potential 
function associated with the Gauss--map 
$\mathcal{G}$, 
given by \[
I:((x_{k},e_{k}))_{k}\mapsto\log\left|\mathcal{G}'\left(\left[\left|x_{1}\right|,\left|x_{2}\right|,\ldots\right]\right)\right|.\]
Also, we require the potential function $J:\Sigma_{G}\rightarrow\R^{2\mathfrak{g}}$ 
given for $\left((x_{k},e_{k})\right)_{k}\in\Sigma_{G}$ by
\[
J\left(\left((x_{k},e_{k})\right)_{k}\right):=\left(\left\langle 
\{ e_{1}\left(i\infty\right),e_{1}(0)\}_{G},f_{1}\right\rangle ,\ldots,\left\langle 
\{ e_{1}\left(i\infty\right),e_{1}(0)\}_{G},f_{2\mathfrak{g}}\right\rangle 
\right),\]
where we will think of $J$ as given by the vector $J=: 
\left(J_{1},\ldots,J_{2\mathfrak{g}}\right)$.\\
Finally, the modular pressure function 
 $P:\R^{2\mathfrak{g}}\times(1/2,\infty)\rightarrow\R$ associated 
 with $J$
is then defined for $t=(t_{1},\ldots,t_{2\mathfrak{g}})\in\R^{2\mathfrak{g}}$ and $\beta\in(1/2,\infty)$
by  (here, $[ \, \, ]$ refers to the cylinder set in $\Sigma_{G}$) \[
P\left(t,\beta\right):=\lim_{n\rightarrow\infty}\frac{1}{n}\log\sum_{\omega\in\Sigma_{G}^{n}}
\exp S_{n}\sup_{x\in[\omega]}\left(\left(t|J(x)\right) -\beta I(x)\right).\]
 Note that $\left(t|J\right) -\beta I$ is acceptable in the
sense of Mauldin/Urba{\'n}ski (\cite[Def. 2.1.4]{MauldinUrbanski:03}),
and this implies that $P$ is well--defined. Also, since $J$ is H\"{o}lder continuous
and bounded, one immediately verifies that $\left(t|J\right)-\beta I$
is summable for each $\beta>1/2$ (for the definition of summablity we refer 
to 
\cite[p. 27]{MauldinUrbanski:03}). In particular, this also gives that  
$P$ is continuous. An argument similar to \cite[paragraph 
6]{KesseboehmerStratmann:04c} then gives that $\lim_{\beta\searrow\frac{1}{2}}P(t,\beta)=\infty$
and $\lim_{\beta\rightarrow\infty}P(t,\beta)=-\infty$. Combining this 
with the continuity of $P$, it follows that there exists a function 
\begin{equation}\label{betafunction} \beta_{G} :\R^{2\mathfrak{g}}\to \left({1}/{2},\infty \right)\end{equation}
such that for each $t\in\R^{2\mathfrak{g}}$ we have $P\left(t,\beta_{G}\left(t\right)\right)=0$. 

We require the following 
facts from the general thermodynamic formalism, which can be found for 
instance  in \cite{MauldinUrbanski:03}. 

\begin{itemize}
\item  For
the potential function $\left(t|J\right)-\beta_{G}(t) I$ there exists a unique ergodic Gibbs measure 
$\mu_{t,\beta_{G}}$  which is positive
on open subsets of $\Sigma_{G}$. In particular, we hence have that 
there exists a constant $Q>1$ such that  for 
each  $\omega \in \Sigma_{G}^{n}$ and $x\in [\omega]$ we have
\begin{eqnarray}\label{gibbs}
Q^{-1} \leq \frac{\mu_{t,\beta_{G}} 
([\omega])}{\exp\left(S_{n}\left(\left(t|J(x)\right) -\beta_{G}(t) 
I(x)\right) - n P(t,\beta_{G}(t))\right) }\leq Q. \end{eqnarray}
For ease of notation, throughout we put $\mu_{t}:=\mu_{t,\beta_{G}}$. 
\item By setting
\[ \partial_{t_i}P\left(t,\beta_{G}(t)\right):= \left.\frac{\partial 
P(t,\beta)}{\partial t_{i}} \right| _{\left(t,\beta_{G}(t)\right)}  
    \;\mbox{ and }\; 
\partial_{\beta}P\left(t,\beta_{G}(t)\right):=  \left.\frac{\partial 
P(t,\beta)}{\partial \beta} \right|_{\left(t,\beta_{G}(t)\right)}, 
 \] we 
have
for all $i \in \{1,\ldots,2 \mathfrak{g}\}$, 
\begin{eqnarray}\label{part} 
    \partial_{t_i}P\left(t,\beta_{G}(t)\right)=
	\int J_i\, d\mu_{t} \;\mbox{ and }\; 
\partial_{\beta}P\left(t,\beta_{G}(t)\right)=
    -\int I\, d\mu_{t}  .
    \end{eqnarray}
 With $\alpha_i(t):=\partial_{t_{i}}\beta_{G}\left(t\right)$, the 
 implicit function theorem then  implies that
 \begin{eqnarray}\label{beta}
\alpha_i(t)=-\frac{ \partial_{t_i}P\left(t,\beta_{G}(t)\right)}{ 
\partial_{\beta}P\left(t,\beta_{G}(t)\right)}=\frac{\int J_{i}\, d\mu_{t}}{\int I\, 
d\mu_{t}}, \textrm{ for all } i \in \{1,\ldots,2 \mathfrak{g}\}.\end{eqnarray}
\end{itemize}
Let us deduce   some further results  crucial for the proof of our 
Main Theorem.
Note that parts of the proof of the following result are inspired by 
a similar argument given in \cite{Lalley:87}.
\begin{prop}\label{hessian} 
The Hessian $\left(\nabla^{2}\beta_{G}\right)(t)$ is strictly positive
definite for all $t\in\R^{2\mathfrak{g}}$. In particular, the function $\beta_{G}:\R^{2\mathfrak{g}}\rightarrow\R$
is strictly convex and the gradient map $\nabla\beta_{G}:\R^{2\mathfrak{g}}\rightarrow\nabla
\beta_{G}\left(\R^{2\mathfrak{g}}\right)$
is a diffeomorphism with a well--defined inverse  $t:\nabla\beta_{G}\left(\R^{2\mathfrak{g}}\right)\rightarrow\R^{2\mathfrak{g}}$. 
\end{prop}
\begin{proof}
As before, let $\mu_{t}:=\mu_{t,\beta_{G}}$ denote the unique Gibbs
measure for the potential function $(t|J)-\beta_{G}(t)I$. Also, for ease of 
exposition,
let $J_{0}:=-I$, as well as $\partial_{0}:=\partial_{\beta}$ and
$\partial_{i}:=\partial_{t_{i}}$, for $i=1,\ldots,2 \mathfrak{g}$.
Since $t\mapsto P(t,\beta_{G}(t))$ defines a constant function, 
its partial derivative
with respect to $t_i$  vanishes for all $i=1,\ldots,2 \mathfrak{g}$. This implies
 \[
\partial_{i}\, P\left(t,\beta_{G}(t)\right)=-\partial_{0}\, P\left(t,\beta_{G}
\left(t\right)\right)\partial_{i}\,\beta_{G}\left(t\right), \textrm{ for 
all  }  
i=1,\ldots,2 \mathfrak{g},\]
 By taking  partial
derivatives  with respect to $t_{j}$ on both sides of this equality we obtain

$\partial_{ij}P(t,\beta_{G}(t)) \partial_{0i}P(t,\beta_{G}(t))\partial_j \beta_{G}(t)$
\[=-(\partial_{0j}P(t,\beta_{G}(t)))+
\partial_{00}P(t,\beta_{G}(t)\partial_j \beta_{G}(t))\partial_i \beta_{G}(t)-\partial_{0}\,
 P\left(t,\beta_{G}\left(t\right)\right)\partial_{ij}\,\beta_{G}
\left(t\right)
\]
Hence, by defining \[
A:=\left(\partial_{ij}\, P\left(t,\beta_{G}\left(t\right)\right)\right)_{i,j=
0,\ldots,2 \mathfrak{g}} \:\textrm{and  } C:=\left(c_{ij}\right)_{i,j},\]
 where \[
 c_{ij}:=\left\{ \begin{array}{lll}
\alpha_{j}(t) & \,\textrm{for} & i=0,\, j=1,\ldots,2 \mathfrak{g}\\
\delta_{ij} & \,\textrm{for} & i,j=1,\ldots,2 \mathfrak{g}\end{array}\right.,\]
we obtain
\[
\left(-\partial_{0}\, P\left(t,\beta_{G}\left(t\right)\right)\right)\left(\partial_{ij}\,\beta_{G}
\left(t\right)\right)_{i,j=1,\ldots,2 \mathfrak{g}}=:B=C^{T}AC.\]
Using (\ref{part}),
it is now sufficient to show that $B$ is strictly positive definite,
or what is equivalent that $A$ is positive definite on the image
$\mathrm{Im}(C)$ of $C$. Here, \[\mathrm{Im}(C):=\left\{ \left(\sum_{i=1}^{2
\mathfrak{g}}\lambda_{i}\alpha_{i},\lambda_{1},\ldots,\lambda_{2\mathfrak{g}}\right):\left(
\lambda_{1},\ldots,\lambda_{2\mathfrak{g}}\right)\in\R^{2\mathfrak{g}}\right\}.\]
For this it is sufficient to show that we have for all $y=(y_{0},y_{1},\ldots,y_{2\mathfrak{g}})\in\textrm{Im}
\left(C\right)\setminus\left\{ 0\right\} $,
\[
y^{T}Ay>0.\]
In order to prove this, note that  by 
\cite[Proposition 2.6.14]{MauldinUrbanski:03} we have \begin{eqnarray*}
\partial_{ij}\, P\left(t,\beta_{G}\left(t\right)\right) & = & \sum_{k=0}^{\infty}\mu_{t}\left((J_{i}-\mu_{t}(J_{i}))(J_{j}\circ\sigma^{k}-\mu_{t}(J_{j}))\right)\\
 & = & \sum_{k=0}^{\infty}\mu_{t}\left((J_{j}-\mu_{t}(J_{j}))(J_{i}\circ\sigma^{k}-
 \mu_{t}(J_{i}))\right)=:\sigma_{t}^{2}\left(J_{i},J_{j}\right).\end{eqnarray*}
 Using this, it follows 
 \begin{eqnarray*}
y^{T}Ay & = & \sum_{i,j=0}^{2\mathfrak{g}}y_{i}y_{j}\sigma_{t}^{2}\left(J_{i},J_{j}\right)=\sum_{i,j=0}^{2\mathfrak{g}}\sigma_{t}^{2}\left(y_{i}J_{i},y_{j}J_{j}\right)=\sigma_{t}^{2}\left(\sum_{i=0}^{2\mathfrak{g}}y_{i}J_{i},\sum_{i=0}^{2\mathfrak{g}}y_{i}J_{i}\right)\\
 & =: & \sigma_{t}^{2}\left(\sum_{i=0}^{2\mathfrak{g}}y_{i}J_{i}\right)\geq 0.\end{eqnarray*}
 Since $\alpha_{i}=\mu_{t}(J_{i})/\mu_{t}(I)$, we have for $y=\left(\sum\lambda_{i}\alpha_{i},\lambda_{1},\ldots,\lambda_{2\mathfrak{g}}\right)\in\mathrm{Im}(C)$,
\[
\mu_{t}\left(\sum_{i=1}^{2\mathfrak{g}}\lambda_{i}J_{i} -
\sum_{i=1}^{2\mathfrak{g}}\lambda_{i}\alpha_{i}I\right)=0.\]
 Let us assume by way of contradiction that $\sigma_{t}^{2}\left(\sum_{i=0}^{2\mathfrak{g}}y_{i}J_{i}\right)=0$.
Note that $\sigma_{t}^{2}\left(\sum_{i=0}^{2\mathfrak{g}}y_{i}J_{i}\right)=0$
if and only if $\sum_{i=0}^{2\mathfrak{g}}y_{i}J_{i}$ is cohomologous to $0$
within the class of bounded H\"{o}lder continuous functions. The 
latter means  that there exists a bounded H\"{o}lder continuous function $u$ on
$\Sigma_{G}$ such that (cf. \cite[Lemma 4.8.8]{MauldinUrbanski:03})
\begin{equation}
\sum_{i=0}^{2\mathfrak{g}}y_{i}J_{i}=u-u\circ\sigma.\label{eq:Cohomolog}\end{equation}
 Hence, it remains to show that (\ref{eq:Cohomolog}) implies $(\lambda_{1},\ldots,\lambda_{2\mathfrak{g}})=0$.
In order to see this, we distinguish the following two cases. First,
if $\sum_{i=1}^{2\mathfrak{g}}\lambda_{i}\alpha_{i}\not=0$ then 
$\sum_{i=1}^{2\mathfrak{g}}\lambda_{i}J_{i}
- \sum_{i=1}^{2\mathfrak{g}}\lambda_{i}\alpha_{i}I$
is an unbounded function (since $I$ is unbounded). Since the right
hand side of (\ref{eq:Cohomolog}) is bounded, we then 
immediately have a contradiction.
Secondly, if $\sum_{i=1}^{2\mathfrak{g}}\lambda_{i}\alpha_{i}=0$ then consider
$F:=\sum_{i=1}^{2\mathfrak{g}}\lambda_{i}J_{i}$. We first investigate how $F$
behaves on elements $\omega:=((x_{k},e_{k}))_{k}\in\Sigma_{G}$ which
are periodic in the second coordinate, that is where there exist $p\in\N$
such that $e_{mp+j}=e_{j}$ for all $m=0,1,\ldots$ and $j=1,\ldots,p$.
In this situation we necessarily have that $S_{p-1}F(\omega)=0$,
since otherwise we would have $\lim_{m\rightarrow\infty}|S_{mp-1}F(\omega)|=\lim_{m\rightarrow\infty}|mS_{p-1}F(\omega)|=\infty$
which contradicts (\ref{eq:Cohomolog}). Therefore, \begin{equation}
S_{p-1}F(\omega)=\sum_{j=1}^{2\mathfrak{g}}\lambda_{j}\bigg\langle\sum_{k=1}^{p}\left\{ e_{k}\left(i\infty\right),e_{k}(0)\right\} _{G},f_{j}\bigg\rangle=0.\label{eq:Cohomolog2}\end{equation}
 Now let $\left\{ \gamma_{1},\ldots,\gamma_{2\mathfrak{g}}\right\} $ be a basis
of $H_{1}\left(M_{G},\mathbb{R}\right)$ consisting of cycles. 
Each $\gamma_{i}$ can be represented by an oriented
closed geodesic in $M_{G}$. The forward directions of these geodesics
correspond to elements  $z_{i}=\left[x_{1},x_{2},\ldots\right] \in 
[0,1]$ which  are periodic in their continued fraction expansion, of period
$2r_{i}$ say (if the period is odd, then replace the geodesic by twice
the geodesic. Since $E_{G}$ is finite, it follows that there
exists $m_{i}\in\N$ and $e_{i,1}\in E_{G}$ such that for $\omega_{i}:=
\left(\left(-x_{1},e_{i,1}\right),\left(x_{2},e_{i,2}\right),\left(-x_{3},e_{i,3}\right),\ldots\right)\in\Sigma_{G}$
we have $e_{i,2m_{i}r_{i}k}=e_{i,1}$, for all $k\in\N$ (note, in this 
step it is vital that the periods were chosen to be even). Hence, $\omega_{i}$
is of period $2m_{i}r_{i}$. This
shows that  the set
\[
\left\{ \sum_{k=1}^{2m_{i}r_{i}}\left\{ 
e_{i,k}\left(i\infty\right),e_{i,k}(0)\right\} 
_{G}:i=1,\ldots,2 \mathfrak{g}\right\} \]
contains a basis of $H_{1}\left(M_{G},\mathbb{R}\right)$.
 The assertion now follows from combining  (\ref{period}) and (\ref{eq:Cohomolog2}). 
\end{proof}
The following immediate corollary   shows that $\beta_{G}$ and its Legendre 
transform $\hat{\beta_{G}}$, given by \[\hat{\beta_{G}}(\alpha):=\inf_{t\in \R^{2\mathfrak{g}}} 
    \left\{\beta_{G}(t)-\left(t|\alpha\right) \right\},\]
    are in fact a Legendre transform pair (cf. \cite{Rockafellar:70}). 
\begin{cor}\label{cor:legendre}
    For the Legendre transform $\hat{\beta_{G}}$ of $\beta_{G}$,  we
    have for each $\alpha \in \nabla\beta_{G}(\R^{2\mathfrak{g}})$, \[
    \hat{\beta_{G}}\left(\alpha\right) =\beta_{G}(t(\alpha)) - 
    \left(t(\alpha)|\alpha\right).\]
\end{cor}
For the final proposition of this section we require the following lemma. 

\begin{lem}\label{product}
For any measure $\mu\in\mathcal{M}\left(\Sigma_{G},\sigma\right)$
we have that the first marginal of $\tilde{\mu}:=\mu\circ\tilde{\pi}^{-1}$
is a shift invariant measure on $\left(\Sigma_{*},\sigma_{*} \right)$.
Furthermore, if $\tilde{\mu}$ can be written as a product measure
$\nu\otimes \mathbb{P}$ on $\Sigma_{*}\times E_{G}$ such that $\nu\left(U\right)>0$
 for all non-empty open subsets $U\subset \Sigma_{*}$, then $\mathbb{P}$ is equal
to the equidistribution on $E_{G}$, that is $\mathbb{P}\left(\left\{ e\right\} \right)=1/\kappa$,
for all $e \in E_{G}$. 
\end{lem}
\begin{rem*}Note that the first part of this lemma in particular shows 
that $\nu$ is a ${\sigma}_*$--invariant measure on $\Sigma_{*}$. It is then an 
immediate consequence of the ergodic theorem and the symmetry
of $E_{G}$, that  the limiting symbol
vanishes almost surely for product
measures of this type. In fact, this special situation occurs for the 
generalized Gauss--measure as  discussed in
\cite{ManinMarcolli:02}, as well as for the Lyapunov spectrum 
arising from continued fraction expansions with bounded entries as 
studied in  \cite{Marcolli:03}. 
\end{rem*}
\begin{proof}[Proof of Lemma \ref{product}]
For the first part let $A\subset \Sigma_{*}$ be some given Borel set. We then
have \begin{eqnarray*}
\tilde{\mu}\left(A\times E_{G}\right) & = & \tilde{\mu}\left(\tilde{\sigma}^{-1}
\left(A\times 
E_{G}\right)\right)=\tilde{\mu}\left(\left(\hat{\sigma}^{-1}A\right)\times E_{G}\right).
\end{eqnarray*}
 For the second part let $p_{e}:=\mathbb{P}\left(\{ e\}\right)$ for $e \in 
 E_{G}$, and define
$p:=\max\left\{ p_{e}: e \in E_{G}\right\} $.  Also, let  $e'$ 
refer to some element of  $E_{G}$ such that $p_{e'}=p$, and define
 for $m\in\N$ 
and $e 
\in E_{G}$, \[
C_{e,e'}^{m}:=\left\{ x\in\left(\Z^{^{\times}}\right)^{\N}:\tau_{x_{1}}\cdots
\tau_{x_{m}}\left(e\right)=e'\right\}  .\]
 For $n$ greater than the maximal word length of 
 the elements in $W$, the 
 $\tilde{\sigma}$--invariance of $\tilde{\mu}$ then
 gives
\begin{eqnarray*}
p & = & p_{e'} = 
\frac{1}{n}\sum_{m=1}^{n}\tilde{\mu}\left(\tilde{\sigma}^{-m}\left(\Sigma_{*}
\times\left\{ e'\right\} \right)\right)=\frac{1}{n}\sum_{m=1}^{n}\tilde{\mu}
\left(\bigcup_{e\in E_{G}}C_{e,e'}^{m}\times\left\{ e\right\} 
\right)\\ &=&\frac{1}{n}\sum_{m=1}^{n}
\sum_{e\in E_{G}}\tilde{\mu}\left(C_{e,e'}^{m}\times\left\{ e\right\} \right)
  = \sum_{e\in 
  E_{G}}\frac{1}{n}\sum_{m=1}^{n}\nu\left(C_{e,e'}^{m}\right)\cdot p_{e}.\end{eqnarray*}
 Now, note that $\sum_{e\in 
  E_{G}}\nu\left(C_{e,e'}^{m}\right)=1$, and therefore  $\sum_{e \in 
 E_{G}}\frac{1}{n}\sum_{m=1}^{n}\nu\left(C_{e,e'}^{m}\right)=1$. 
Combining this with the fact that by assumption we have 
$\frac{1}{n}\sum_{m=1}^{n}\nu\left(C_{e,e'}^{m}\right)>0$, the 
calculation above implies that
 $p_{e}=p$ for each $e \in E_{G}$. 
\end{proof}

\begin{prop}
\label{lem:vanishing}
We have  $\partial_{t_{i}}P\left(0,\beta_{G}\left(0\right)\right)=0$, 
for all $i\in\left\{ 1,\ldots,2 \mathfrak{g}\right\} $. 
\end{prop}
\begin{proof}
 Since $\mu_{0}$ is the unique ergodic Gibbs measure
for the potential function $-I$, it follows from Lemma \ref{product} 
that  the
pull--back $\mu_{0}\circ\tilde{\pi}^{-1}$ of $\mu_{0}$ to  
$\tilde{\Sigma}_{G}$ 
can be written as a product measure $\nu \otimes\mathbb{P}$, where
 $\mathbb{P}$ refers to the equidistribution
on $E_{G}$. Next note  that with $S$ referring to the elliptic 
generator of $\Gamma$ of order $2$,  we have that  $\left\{ eS:e\in 
E_{G}\right\} $
is a set of representatives of $G \backslash \Gamma$. Using this, it follows  \begin{eqnarray*}
\partial_{t_{i}}P\left(0,\beta_{G}\left(0\right)\right) & = & \int\left\langle 
\left\{ e\left(i\infty\right),e\left(0\right)\right\} _{G},f_{i}\right\rangle \, d\mu_{0}\\
&=&\frac{1}{\kappa}\sum_{e\in E_{G}}\left\langle \left\{ e\left(i\infty\right),e\left(0\right)\right\} _{G},f_{i}\right\rangle \\
 & = & \frac{1}{\kappa}\sum_{e\in E_{G}}\left\langle \left\{ 
 eS\left(i\infty\right),eS\left(0\right)\right\} _{G},f_{i}\right\rangle \\
 &=&\frac{1}{\kappa}\sum_{e\in E_{G}}\left\langle \left\{ e\left(0\right),e\left(i\infty\right)\right\} _{G},f_{i}\right\rangle \\
 & = & -\frac{1}{\kappa}\sum_{e\in E_{G}}\left\langle \left\{ e\left(i\infty\right),e\left(0\right)\right\} _{G},f_{i}\right\rangle =-\partial_{t_{i}}P\left(0,\beta_{G}\left(0\right)\right).\end{eqnarray*}
This implies that $\partial_{t_{i}}P\left(0,\beta_{G}\left(0\right)\right)=0$. 
\end{proof}
\begin{cor}\label{cor:beta}
The function $\beta_{G}:\R^{2\mathfrak{g}}\rightarrow\R$ has a unique minimum 
at $0\in\R^{2\mathfrak{g}}$, and  $\beta_{G}(0)=1$.
\end{cor}
\begin{proof}
Combining (\ref{beta}) and Proposition \ref{lem:vanishing},  it follows that $\nabla\beta_{G}\left(0\right)=0$.
Also, by Proposition \ref{hessian} we have that $\beta_{G}$ is strictly 
convex. Combining these two observations,  
it follows that $\beta_{G}$ has a  unique minimum at zero. Finally, 
note that $P(0,1)=0$ (see \cite{KesseboehmerStratmann:04c}),
which immediately implies that $\beta(0)=1$.
\end{proof}

\section{Proof of Main Theorem}
For $\alpha\in \R^{2\mathfrak{g}}$
and $e_{1}\in E_{G}$,  consider the set\[
\mathcal{F}_{\alpha}(e_{1}):=\left\{ 
x\in \mathcal{I} : ((x_{k},e_{k}))_{k} \in \Sigma_{G} \, \textrm{ such 
that } \, \lim_{n\to\infty}\frac{S_n J
(((x_{k},e_{k}))_{k})}{S_n I(((x_{k},e_{k}))_{k})}=\alpha \right\} ,\]
where as before, we have set
 $x:=-\sign\left(x_{1}\right) \left[\left|x_{1}\right|,\left|x_{2}\right|,\ldots \right]$.
One immediately verifies that 
$\left(q_n(|x|)\right)^{2} \asymp \exp \left( S_n 
I\left(((x_{k},e_{k}))_{k}\right)\right)$, and hence 
 by 
 Proposition \ref{partial} we have for each 
 $((x_{k},e_{k}))_{k} \in \Sigma_{G}$ and $\alpha \in \R ^{2\mathfrak{g}}$,
 \[
  \left(\langle\tilde{\ell}_G(((x_{k},e_{k}))_{k}),f_{1}\rangle,\ldots ,
  \langle\tilde{\ell}_G(((x_{k},e_{k}))_{k}),f_{2\mathfrak{g}}\rangle
  \right)=\alpha \iff \lim _{n\to \infty} 
  \frac{S_n J (((x_{k},e_{k}))_{k})}{S_n I(((x_{k},e_{k}))_{k})}=\alpha.\] 
  This shows   that
 $ \mathcal{F}_{\alpha}(e) = \tilde{\mathcal{F}}_{\alpha}(e)$, for 
 all $e \in E_{G}$. Hence, by Lemma \ref{eql}, it is now sufficient  to compute 
$\dim_{H}\left(\mathcal{F}_{\alpha}(e)\right)$ for some fixed $e \in E_{G}$.
 For this, we define for  $y =[y_{1},y_{2},\ldots]  \in 
 \mathcal{I}_{+1}$ and with $B(y,r)$ referring to the interval 
 centred at $y$ of radius $r>0$,  
    \[n_{r}(y):= \min \{n : [y_{1},\ldots,y_{n}] \subset B(y,r) 
\},\;\;
m_{r}(y):= \max \{n : [y_{1},\ldots,y_{n}] \supset B(y,r) 
\}. \]
Note that we clearly have that $|n_r(y)-m_r(y)|$ is uniformly bounded from above.  
Using this and the Gibbs property of the pull--back 
$\overline{\mu}_{t(\alpha)}$ of $\mu_{t(\alpha)}$ to 
$\overline{\Sigma}_{G}$ (see (\ref{gibbs})), we then have  
for each $\alpha\in\nabla\beta_{G}\left(
\R^{2\mathfrak{g}}\right) $ and for $\overline{\mu}_{t(\alpha)}$--almost every 
$(e(x),e) \in \overline{\Sigma}_{G}$, where as before  $x =-\sign\left(x_{1}\right) 
\left[\left|x_{1}\right|,\left|x_{2}\right|,\ldots \right] \in \mathcal{F}_{\alpha}(e)$
and $((x_{k},e_{k}))_{k}$ refers to the corresponding element in $\Sigma_{G}$,
\begin{eqnarray*}   & & \,  \hspace{-12mm} \lim_{r\to 
0}   \frac{\log\overline{\mu}_{t(\alpha)}(B(e(x),r) \times \{e\})}{\log 
r}\\
&=&
\lim_{r\to 0} \frac{ 
\left(t(\alpha)|S_{n_{r}(|x|)}J\left(\left((x_{k},e_{k})\right)_{k}\right)\right)
-\beta(t(\alpha)) S_{n_{r}(|x|)} 
I\left(((x_{k},e_{k}))_{k}\right) }{-S_{n_{r}(|x|)}
I\left(((x_{k},e_{k}))_{k}\right)} \\
&=&   
\beta_{G}(t(\alpha)) - (t(\alpha)|\alpha)= \hat{\beta_{G}}(\alpha) ,
\end{eqnarray*}
where the last equality follows from Corollary \ref{cor:legendre}. 
Note that by combining (\ref{beta}) and the ergodicity of $\mu_t$,  we 
have that
$\overline{\mu}_{t(\alpha)}(e\left(\mathcal{F}_{\alpha}(e) \right)\times \{e\})/
\overline{\mu}_{t(\alpha)}
(e\left(\mathcal{I}\right) \times \{e\})=1$. 
Therefore, we can apply the mass distribution principle (cf. 
\cite{Falconer:90}) which gives 
\[
\dim_{H}\left(\mathcal{F}_{\alpha}(e)\right)=\hat{\beta_{G}}(\alpha).\]
The remaining assertions of the Main Theorem are 
obtained as follows.  
Since $\mathcal{F}_{\alpha}$ can be written as a countable union of conformal 
images of the sets $ \mathcal{F}_{\alpha}(e)$ for $e\in E$, 
  we have that $\dim(\mathcal{F}_{\alpha})= 
\dim_H(\mathcal{F}_{\alpha}(e))$. This shows that 
$\dim(\mathcal{F}_{\alpha}) = \hat{\beta}_{G}(\alpha)$ for all 
$\alpha\in\nabla\beta_{G}(\R^{2\mathfrak{g}})$, and hence gives the 
assertion in (\ref{mainbeta}). The facts $\beta_{G}(0)=1$  and that $\beta_{G}$ has a 
unique minimum at $0\in \R^{2\mathfrak{g}}$ have already been obtained in Corollary 
\ref{cor:beta}. Likewise, 
the analytic properties of $\beta_G$ stated 
in the Main Theorem can be deduced from Proposition \ref{hessian}. 
Finally, for the assertion that the dimension spectrum is almost complete 
we refer to \cite[Theorem 1.2]{KessUrbanski:06} (note that
the results  
in \cite{KessUrbanski:06} are applicable since $\beta_{G}$ is strictly convex).

\def\cprime{\('\)}

\end{document}